\documentclass[letterpaper,12pt]{article}

\usepackage{color}
\usepackage{amssymb}
\usepackage{amsmath}
\usepackage{amsthm}
\usepackage{amscd}
\usepackage{mathtools}
\usepackage[all]{xy}
\usepackage{tikz-cd}

\usepackage{hyperref}
\hypersetup{colorlinks=true,linkcolor=blue,anchorcolor=blue,citecolor=blue}
\usepackage{cleveref}

\theoremstyle{plain}
\newtheorem{theo}{Theorem}[section]
\newtheorem{prop}[theo]{Proposition}
\newtheorem{lem}[theo]{Lemma}
\newtheorem{cor}[theo]{Corollary}
\theoremstyle{definition}
\newtheorem{defi}[theo]{Definition}
\newtheorem{rem}[theo]{Remark}

\renewcommand \deg {{\mathrm{deg}}}
\newcommand \id {\mathrm{id}}
\newcommand \Br {{\mathrm{Br}}}

\newcommand \si {{\sigma}}
\newcommand \Ga {{\Gamma}}
\newcommand \R {{\mathbb{R}}}
\newcommand \Pic {\mathrm{{Pic}}}
\newcommand \Div {\mathrm{{Div}}}

\newcommand \Gal {{\mathrm{Gal}}}
\newcommand \Ker {{\mathrm{Ker}}}
\renewcommand \Im {\mathrm{{Im}}}
\newcommand \Frob{{\mathrm{Frob}}}

\newcommand \A{{\mathbb A}}
\renewcommand \P{{\mathbb P}}

\newcommand \Hom {\mathrm{{Hom}}}

\newcommand \Cl{\mathrm{Cl}}

\newcommand\ov{\overline}

\newcommand \Z {{\mathbb Z}}
\newcommand \Q {{\mathbb Q}}
\newcommand \F {{\mathbb F}}

\newcommand\G{{\mathbb G}}

\newcommand\C{{\mathbb C}}

\newcommand\lra{\longrightarrow}

\renewcommand\H{\mathrm{H}}

\newcommand\sA{\mathcal{A}}

\newcommand\Tor{\mathrm{Tor}}

\newcommand\NS{\mathrm{NS\,}}

\renewcommand\si{\sigma}

\renewcommand\Ga{\Gamma}

\newcommand\e{\varepsilon}
\newcommand\et{{\mathrm{\acute et}}}

\newcommand \g{{\mathfrak g}}

\newcommand\fH{{\mathfrak H}}

\newcommand{\bthe}{\begin{theo}}
\newcommand{\ble}{\begin{lem}}
\newcommand{\bpr}{\begin{prop}}
\newcommand{\bco}{\begin{cor}}
\newcommand{\bde}{\begin{defi}}
\newcommand{\brem}{\begin{rem}}
\newcommand{\bprf}{\begin{proof}}
\newcommand{\ethe}{\end{theo}}
\newcommand{\ele}{\end{lem}}
\newcommand{\epr}{\end{prop}}
\newcommand{\eco}{\end{cor}}
\newcommand{\ede}{\end{defi}}
\newcommand{\erem}{\end{rem}}
\newcommand{\eprf}{\end{proof}}

\numberwithin{equation}{section}

\title{Surfaces defined by pairs of polynomials}
\author{Dami\'an Gvirtz-Chen and Alexei N.~Skorobogatov}
\begin{document}
\maketitle

\begin{abstract}
We compute the Brauer group of surfaces defined by equating
two bilinear forms of the same degree, assuming these forms are, in an explicit sense, sufficiently general. Our method uses a topological deformation argument
and does not require full knowledge of the algebraic or transcendental cycles.
We obtain a criterion for the triviality of the transcendental Brauer group of an isotrivial variety, 
which we use to prove that
the Brauer group of the generic diagonal surface of arbitrary degree is trivial.
\end{abstract}


\section{Introduction}

Let $k$ be a field of characteristic zero. Consider a smooth surface
$X\subset\P^3_k$ given by
\begin{equation}
F(x_0,x_1)=G(x_2,x_3), \label{F=G}
\end{equation}
where $F$ and $G$ are homogeneous forms of degree $d$ without multiple roots. 
The surface $X$ is geometrically rational when $d\leq 3$, a K3 surface when $d=4$, and a surface
of general type when $d\geq 5$.
The surface $X$ is covered by the product of the plane curves  $y^d=F(x_0,x_1)$ and $z^d=G(x_2,x_3)$, which implies that the Brauer group $\Br(X)$
is finite modulo the image of $\Br(k)$ when $k$ is finitely generated over $\Q$, see \cite[Cor.~1.8]{GS22}.

Let $\bar k$ be an algebraic closure of $k$, let 
$\Ga_k=\Gal(\bar k/k)$, and let $X_{\bar k}=X\times_k\bar k$.
The Brauer group $\Br(X)$ has a natural filtration
\[
\Br_0(X)=\Im[\Br(k)\to\Br(X)], \quad \Br_1(X)=\Ker[\Br(X)\to \Br(X_{\bar k})].
\]
We denote by $d_{1,1}\colon\H^1(k,\Pic(X_{\bar k})) \to \H^3(k,\bar k^\times)$
the differential in the Hochschild--Serre spectral sequence 
\begin{equation}
E^{p,q}_2\colon \ \H^p(k,\H^q(X_{\bar k}, \G_m))\Rightarrow \H^{p+q}(X,\G_m). \label{ss}
\end{equation}
\sloppy As is well known, the spectral sequence (\ref{ss}) implies that 
$\Br_1(X)/\Br_0(X) \cong \Ker(d_{1,1})$.
Determining $\Br_1(X)/\Br_0(X)$
(respectively, $\Br(X)/\Br_1(X)$) usually requires explicit knowledge of the Galois module structure
of $\Pic(X_{\bar k})$
(respectively, of $\Pic(X_{\bar k})$ and $T(X_{\bar k})\otimes\Q/\Z$,
where $T(X_{\bar k})$ is the group of the transcendental cycles).

In our first result we completely determine $\Br_1(X)/\Br_0(X)$ for the surfaces given by 
(\ref{F=G}) when $k$ is a Hilbertian field (for example, a field finitely generated over $\Q$ \cite[\S 9.6]{Serre})
and the forms $F$ and $G$ are sufficiently general. Let us emphasise that for $d\geq 4$ these surfaces
do not lie in a single geometric isomorphism class and we do not know the explicit structure of
$\Pic(X_{\bar k})$ as an abelian group and a Galois module.

Assume that $F(x_0,x_1)$ and $G(x_2,x_3)$ are {\em irreducible} over $k$.
Then $f(t)=F(1,t)$ and $g(t)=G(1,t)$ are irreducible polynomials in $k[t]$ of degree $d$.
Let $k_f$ and $k_g$ be the splitting fields and let $G_f=\Gal(k_f/k)$ and $G_g=\Gal(k_g/k)$
be the Galois groups of $f(t)$ and $g(t)$, respectively.
Assume that $k_f$ and $k_g$ are {\em linearly disjoint over $k$}, that is, $k_f\cap k_g=k$. 
Fix a root $\xi\in\bar k$ of $f(t)=0$ and consider the subgroup $S_f=\Gal(k_f/k(\xi))\subset G_f$ 
and its normal closure $N_f\subset G_f$. Define $S_g$ and $N_g$ similarly.

Let $G=G_f\times G_g$. 
For a $G$-module $M$ we denote by $\H^i(G,M)_\mathrm{prim}$ the kernel of the restriction map
$\H^i(G,M)\to\H^i(G_f,M)\oplus \H^i(G_g,M)$. 
The K\"unneth formula for group homology and the universal coefficients formula
give rise to an isomorphism \[\H^2(G,\Z/d)_\mathrm{prim}\cong
\Hom(G_f^\mathrm{ab}\otimes G_g^\mathrm{ab},\Z/d),\] see (5.36) on p.~154 of \cite{CTS21}.

Let $\varphi$ be a bilinear form
$(G_f/N_f)^\mathrm{ab}\times (G_g/N_g)^\mathrm{ab}\to\Z/d$.
Using the canonical surjective function $G_f/S_f\to (G_f/N_f)^\mathrm{ab}$ and a similar map for $g$,
we obtain from $\varphi$ a function
$\widetilde\varphi\colon (G_f/S_f)\times (G_g/S_g)\to\Z/d$.
Define 
\[
\Pi_{f,g}\subset \Hom((G_f/N_f)^\mathrm{ab}\otimes (G_g/N_g)^\mathrm{ab},\Z/d)
\]
to be the subgroup consisting of the bilinear forms $\varphi$ satisfying
\[
\sum_{y\in G_g/S_g}\widetilde\varphi(\alpha,y)=
\sum_{x\in G_f/S_f} \widetilde\varphi(x,\beta)=0
\]
for any $\alpha\in G_f/S_f$ and any $\beta\in G_g/S_g$.
For $a\in k^\times$ let $\Delta_a$ be the composition 
\[
\Pi_{f,g} \hookrightarrow\Hom((G_f/N_f)^\mathrm{ab}\otimes (G_g/N_g)^\mathrm{ab},\Z/d)\hookrightarrow
\Hom(G_f^\mathrm{ab}\otimes G_g^\mathrm{ab},\Z/d)
\]
\[
\tilde\lra\H^2(G,\Z/d)_\mathrm{prim}\hookrightarrow\H^2(G,\Z/d)\stackrel{\mathrm{inf}}\lra
\H^2(k,\Z/d)\stackrel{\cup[a]}\lra\H^3(k,\mu_d)\to\H^3(k,\bar k^\times),
\]
where $[a]$ denotes the class of $a$ in $\H^1(k,\mu_d)=k^\times/k^{\times d}$.
If the orders of $(G_f/N_f)^\mathrm{ab}$ and $(G_g/N_g)^\mathrm{ab}$ are coprime, then $\Pi_{f,g}=0$.
For more cases when $\Pi_{f,g}=0$ see Proposition \ref{Pi} below.

Our main result for the case of Hilbertian fields is the following theorem.

\bthe \label{general}
Let $k$ be a Hilbertian field of characteristic zero, for example, a field finitely generated over $\Q$.
Let $F(x_0,x_1)$ and $G(x_2,x_3)$ be irreducible binary forms of degree $d$
over $k$ whose splitting fields $k_f$ and $k_g$ are linearly disjoint over $k$. 
Assume that the Jacobians of the plane curves
$y^d=F(x_0,x_1)$ and $z^d=G(x_2,x_3)$ are not isogenous over $\bar k$.
Then the surface $X\subset\P^3_k$ given by $F(x_0,x_1)=G(x_2,x_3)$ satisfies the
following properties.

\textrm{(i)} $\H^1(k,\Pic(X_{\bar k}))\cong \Pi_{f,g}$. 

\textrm{(ii)} Let $K$ be the compositum of $k(\mu_d)$, $k_f$, $k_g$.
Then the differential $d_{1,1}$ can be identified with $\Delta_a$,
where $a=G(u)/F(u)\in k^\times$ with $u\in\P^3_k(k)$ is such that
the restriction of $[a]$ to $\H^1(K,\mu_d)$ has order $d$.
\ethe
\brem{\label{Hilb}
We can always find a $k$-point $u\in\P^3_k(k)$ such that the restriction
of the ratio $[G(u)/F(u)]\in \H^1(k,\mu_d)$ to $\H^1(K,\mu_d)$ has order $d$. Indeed,
by \cite[\S 9.2, Proposition 1]{Serre}
there is a thin set $\Omega\subset\P^3_k(K)$ such that for all $u\in\P^3_k(K)\setminus\Omega$ 
the polynomial $t^d-G(u)/F(u)$ is irreducible over $K$.
Since $k(\mu_d)\subset K$, this implies that the class of $G(u)/F(u)$ in 
$\H^1(K,\mu_d)$ has order $d$.
By \cite[\S 9.4, Proposition]{Serre} the set $\P^3_k(k)\cap\Omega$ is thin over $k$.
Since $k$ is Hilbertian, the complement to a thin subset of $\P^3_k(k)$ is Zariski dense, so
there are infinitely many $k$-points $u\in\P^3_k(k)$ with the required property.
}\erem

The assumption on the Jacobians in Theorem \ref{general} is not essential and can be relaxed as stated in Theorem \ref{f=g}.

\bpr \label{Pi}
Let $f(t)$ and $g(t)$ be irreducible polynomials in $k[t]$ of degree $d$ such that their
splitting fields are linearly disjoint over $k$. 
If at least one of $f$ and $g$ has primitive Galois group,
then $\Pi_{f,g}=0$, except when $d$ is an odd prime
and both polynomials have cyclic Galois groups.
For any $d$, if $G_f\simeq G_g\simeq\Z/d$, then $\Pi_{f,g}\cong\Z/d$ if $d$ is odd, and 
$\Pi_{f,g}\cong\Z/(d/2)$ if $d$ is even.
\epr

Under the assumptions of Theorem \ref{general}, we deduce from Proposition \ref{Pi}
that $\Br_1(X)=\Br_0(X)$ if at least one of 
$G_f$ and $G_g$ is the symmetric group $S_d$ for $d\geq 3$, or the alternating group $A_d$ for
$d\geq 4$. 
The same conclusion holds if $d$ is prime and at least one of $G_f$ and $G_g$ is not isomorphic to $\Z/d$, hence singling out the cyclic case as special.

For certain `isotrivial' surfaces over a function field we can prove the triviality of the
transcendental Brauer group. See Theorem \ref{friday}, from which our next result follows.

\bthe \label{n}
Let $k_0$ be a field of characteristic zero, and let $f$ and $g$ be separable
polynomials in $k_0[t]$ of degree $d$. Let $k=k_0(a)$
be a purely transcendental extension of $k_0$ of transcendence degree $1$, and
let $X\subset \P^3_k$ be the surface $(\ref{X})$.
Then the map $\Br(X)\to\Br(X_{\bar k})$ is zero. 
\ethe

\subsubsection*{Diagonal surfaces}

Theorem \ref{general} is deduced from Theorem \ref{f=g}. The latter can also
be applied to diagonal surfaces, that is, the surfaces 
$X\subset\P^3_k$ given by
\begin{equation}
x_0^d+a_1x_1^d+a_2x_2^d+a_3x_3^d=0,\label{eqX}
\end{equation}
where $a_1,a_2,a_3\in k^\times$.  Such surfaces can be regarded as twisted forms of
the Fermat surface $F\subset\P^3_\Q$:
\begin{equation}
x_0^d+x_1^d+x_2^d+x_3^d=0.\label{Fermat}
\end{equation}
Let $L$ be the field of definition of $\Pic(F_{\ov\Q})$, that is, the smallest extension of $\Q$
such that $\Gal(\ov\Q/L)$ acts trivially on $\Pic(F_{\ov\Q})$. In \S\ref{Pic}
we determine the field $L$ for any $d$ via a Größencharacter computation, thereby completing previous work of Shioda and Aoki \cite{Sh81, Aok, ASh}.
Let us denote by $k_L$ the compositum $kL$, and by $\langle a_1,a_2,a_3\rangle$
the subgroup of $k_L^\times/(k_L^\times)^d$ generated by the classes of $a_1,a_2,a_3$.
The following condition describes what we mean by `general' diagonal surfaces:
\[
\langle a_1,a_2,a_3\rangle\simeq(\Z/d)^3. \eqno{(\star)}
\]
Fixing a primitive $d$-th root of unity we have a map
\[
\H^3(k,\mu_d^{\otimes 3})\to \H^3(k(\mu_{d}),\mu_d^{\otimes 3})
\to \H^3(k(\mu_{d}),\bar k^\times),
\]
where the first arrow is restriction. Consider the following condition:
\[
\text{the order of the image of $[-a_1]\cup[-a_2]\cup[-a_3]$ in 
$\H^3(k(\mu_{d}),\bar k^\times)$ is $d$.}\eqno{(\star\star)}
\]

Let $\Q(\mu_{d})^{[2]}\subset \Q(\mu_{d})$ be the invariant subfield of 
the 2-torsion subgroup $(\Z/d)^\times[2]$ of $(\Z/d)^\times\cong \Gal(\Q(\mu_{d})/\Q)$. 

Write $\chi\colon\Ga_k\to (\Z/d)^\times$ for the cyclotomic character of degree $d$.

We have $2(\Z/d)\simeq\Z/d$ when $d$ is odd and $2(\Z/d)\simeq\Z/(d/2)$ when $d$ is even. 

\bthe\label{main}
Let $X\subset \P^3_k$ be the diagonal surface $(\ref{eqX})$
satisfying condition $(\star)$.
Then the following properties hold.

\textrm{(i)} 
The group $\H^1(k,\Pic(X_{\bar k}))$ is isomorphic to the largest subgroup of $2(\Z/d)$ on which
$\chi^2(\Ga_k)$ acts by $1$. In particular, if $\Q(\mu_{d})^{[2]}\subset k$, then 
$\H^1(k,\Pic(X_{\bar k}))\simeq 2(\Z/d)$.

\textrm{(ii)} If $\Q(\mu_{d})\subset k$ and 
$r|d$ is the order of the image of $[-a_1]\cup[-a_2]\cup[-a_3]$ in 
$\H^3(k,\bar k^\times)$, then $\Br_1(X)/\Br_0(X)\cong \H^1(k,\Pic(X_{\bar k}))[d/r]$.
For any $k$, if condition $(\star\star)$ holds, for example when $k=k_0(a_1,a_2,a_3)$
is a purely transcendental extension of transcendence degree $3$ of a subfield $k_0$,
then $\Br_1(X)=\Br_0(X)$.

\textrm{(iii)} The natural map $\Br(k)\to\Br(X)$ is injective.
\ethe

The result of Theorem \ref{main} (i) for $d=3$ is due to J.-L.~Colliot-Th\'el\`ene, 
D.~Kanevsky and J.-J.~Sansuc when $k$ contains $\Q(\mu_3)$ 
\cite[Prop.~1]{CTKS}. This was extended to arbitrary fields of characteristic zero in
\cite[Prop.~6.1]{PT}, see also \cite[Lemma 2.1]{BBL}. For $d=4$ it is due
to M.~Bright \cite{Bri02, Bri11}. 

The condition $(\star\star)$ holds when $k=k_0(a_1,a_2,a_3)$
is a purely transcendental extension of $k_0$ of transcendence degree $3$. 
In this case Theorem \ref{main} (ii)
was proved by T.~Uematsu for $d=3$ \cite[Thm.~1.2]{U} and by T.~Santens for $d=4$
\cite[Thm.~6.2]{San} using \cite{U2} and \cite{GS22}. Our method works for all degrees $d$.

Condition $(\star\star)$ does not hold when $k$ is a number field, 
because in this case $\H^3(k,\bar k^\times)=0$. For number fields,
Theorem \ref{main} takes the following form.

\bco \label{co2}
Let $k$ be a number field, and let
$X\subset \P^3_k$ be the diagonal surface $(\ref{eqX})$
satisfying condition $(\star)$. The map $\Br(k)\to\Br(X)$ is injective.
The group $\Br_1(X)/\Br(k)$ is isomorphic to the largest subgroup of $2(\Z/d)$ on which
$\chi^2(\Ga_k)$ acts by $1$.
If $k$ contains $\Q(\mu_d)^{[2]}$, then $\Br_1(X)/\Br(k)\simeq 2(\Z/d)$.
\eco

Since $\Q(\mu_{d})^{[2]}=\Q$ precisely when
$d\in\{1, 2, 3, 4, 6, 8,12,24\}$,
 for any degree $d$ from this list and for any number field $k$ we have
$\Br_1(X)/\Br(k)\simeq 2(\Z/d)$. It is easy to show that if $k=\Q$, then 
$\Br_1(X)/\Br(\Q)$ is isomorphic to 
$\Z/2^a\times \Z/3^b$, where $a=2$ if $8$ divides $d$, $a=1$ if 
$4$ simply divides $d$, and $a=0$ if $d$ is not a multiple of $4$,
and $b=1$ if $3$ divides $d$ and $b=0$ otherwise. In particular,
if $d$ is not divisible by 3 or 4, then $\Br_1(X)=\Br(\Q)$.

Comparing Corollary \ref{co2} to Theorem \ref{main} (ii) 
we see that in the family of diagonal surfaces (\ref{eqX}) over $\A^3_k$
`uniform generators' of the algebraic Brauer group do not exist.

Finally, we show the triviality of the Brauer group of the \emph{generic} diagonal surface.

\bthe \label{co1}
Let $k_0$ be a field of characteristic zero and let $k=k_0(a_1,a_2,a_3)$
be a purely transcendental extension of $k_0$ of transcendence degree $3$.
Let $X\subset \P^3_k$ be the diagonal surface $(\ref{eqX})$.
Then the map $\Br(k)\to\Br(X)$ is an isomorphism.
\ethe

In this case both conditions $(\star)$ and $(\star\star)$ hold.
This can be compared to `less generic' diagonal surfaces, where these conditions 
fail, see the examples at the end of \S\ref{3}. 

\subsubsection*{Outline of paper}
In \S\ref{1} we discuss the subgroup
$\Lambda\subset \Pic(X_{\bar k})$ generated by the classes of $d^2$ lines contained in
$F(x_0,x_1)=G(x_2,x_3)=0$. The main result of this section is
Theorem \ref{Lambda} which gives a condition for the natural map 
$\H^1(k,\Lambda)\to\H^1(k,\Pic(X_{\bar k}))$ to be an isomorphism. The proof strategy is a topological deformation argument to the diagonal case, where we build on work of Shioda and Degtyarev.
The $\Ga_k$-module $\Lambda$ has a simple explicit description.
Moreover, the field of definition of $\Lambda$ is the compositum of
$k_f$ and $k_g$. This is why we do not need to know the group
$\Pic(X_{\bar k})$ and the action of $\Ga_k$ on it.

The results of \S \ref{1} are translated into explicit algebraic descriptions in \S \ref{diff}.

Theorem \ref{general} is a consequence of Theorem \ref{f=g}, which is our main technical result.
In~\S\ref{proofs} we prove Proposition \ref{Pi}, Theorem \ref{f=g},
and then deduce Theorems \ref{general} and~\ref{main}. 

In \S\ref{3} we prove the triviality of  the transcendental Brauer group
for some `isotrivial' surfaces, If $k$ is finitely generated over $\Q$, then $k$ is Hilbertian.
see Theorem \ref{friday}. Then we deduce Theorems \ref{n} and \ref{co1}.

Finally, \S\ref{Pic} contains results about the field of definition of 
the geometric Picard group of the surfaces that we study, including
a complete (though complicated) description of $L$.

\medskip

The authors thank ICMS Edinburgh and the Bernoulli Centre at EPFL for their hospitality. The authors are very grateful to the referee for the careful reading of the paper.
This work was partially supported by UKRI award UKRI094.

\section{Reduction to $d^2$ lines} \label{1}
We begin by rewriting equation (\ref{F=G}) using monic polynomials in one variable. 
In this section $f(t)$ and $g(t)$ are separable {\em monic} polynomials of degree $d$ with coefficients
in a field $k$ of characteristic zero, and $X\subset\P^3_k$ is the surface given by 
\begin{equation}
f(x_1/x_0)x_0^d=ag(x_3/x_2)x_2^d, \label{X}
\end{equation}
where $a\in k^\times$. It is immediate to check that $X$ is smooth, and hence geometrically integral. 

Define $Z\subset X$ as the union of closed subsets $x_0=0$ and  
$f(x_1/x_0)x_0^d=0$. 
On $X_{\bar k}$, the zero set of $x_0$ 
is a smooth irreducible curve, while the zero set of $f(x_1/x_0)x_0^d$ is the union of $d^2$
lines. Let $U=X\setminus Z$. 

Consider the subgroup $\Lambda\subset\Pic(X_{\bar k})$ generated by the 
classes of $d^2$ lines contained in $Z_{\bar k}$. 
The following lemma gives an explicit description of this Galois module.
In particular, it shows that the Galois module $\Lambda$
does not depend on $a$. 

Write $V_f$ for the closed subscheme of the affine line $\A^1_k$ given by $f(t)=0$, and similarly
for $V_g\subset\A^1_k$. 

For a 0-dimensional $k$-scheme $S$ of finite type we denote by $\Z[S]$ the free abelian group
whose generators bijectively correspond to the $\bar k$-points of $S$. Then $\Z[S]$
is a permutation $\Ga_k$-module.  Let $\si_S\in \Z[S]$ be the sum of the generators. 
We write $\si_f=\si_{V_f}$ and $\si_g=\si_{V_g}$.

\ble \label{otto}
We have a commutative diagram of $\Ga_k$-modules with exact rows
\begin{equation}
\begin{split}
\xymatrix{
0\ar[r]&\bar k^\times\ar[r]& \bar k(X)^\times\ar[r] &
\Div(X_{\bar k})\ar[r]&\Pic(X_{\bar k})\ar[r]& 0\\
0\ar[r]&\bar k^\times\ar[u]^{\mathrm{id}}_\cong\ar[r]& \bar k[U]^\times\ar[r] \ar@{_{(}->}[u]&
\Div_{Z_{\bar k}}(X_{\bar k})\ar@{_{(}->}[u]\ar[r]&\Lambda\ar@{_{(}->}[u]\ar[r]& 0\\
0\ar[r]& \Z\ar[r]^{p_1 \ \ \ \ \ \ \  \ }\ar[u]^a& \Z[V_f]\oplus\Z[V_g]\ar[r]^{p_2 \ \ \ \ \ }\ar[u]& 
(\Z[V_f]\otimes \Z[V_g])\oplus\Z h\ar[r]\ar[u]_{\mathrm{id}}^\cong& \Lambda\ar[r]\ar[u]_{\mathrm{id}}^\cong& 0
}
\label{duo}
\end{split}
\end{equation}
The arrows from the middle row to the top row are the natural embeddings. The remaining
arrows are as follows:
\begin{itemize}
\item the vertical map with label $a$ sends $1\in\Z$ to $a\in\bar k^\times$;
\item $p_1(1)=\si_f-\si_g$;
\item $h\in\Z h$ is the hyperplane section $x_0=0$;
\item generators $\xi\in\Z[V_f]$, where $f(\xi)=0$, are the rational functions
$(x_1-\xi x_0)/x_0$;
\item generators $\rho\in\Z[V_g]$, where $g(\rho)=0$, are the rational functions
$(x_3-\rho x_2)/x_0$;
\item generators $\xi\otimes\rho$ of $\Z[V_f]\otimes \Z[V_g]$ are the lines 
$x_1-\xi x_0=x_3-\rho x_2=0$;
\item $p_2$ sends each $\xi\in \Z[V_f]$ to $(\xi\otimes\si_g) -h$ and
each $\rho\in \Z[V_g]$ to $(\si_f\otimes\rho)-h$;
\item unlabeled horizontal arrows are the natural maps.
\end{itemize}
\ele
\bprf The commutativity of the diagram is easily checked. The exactness of the top and
middle sequences is clear. It remains to prove the exactness of the bottom sequence, which
is non-trivial only at $(\Z[V_f]\otimes \Z[V_g])\oplus\Z h$. We need to check
that every relation among the classes
of the $d^2$  lines in $\Lambda$ is a linear combination of the relations given by
the rational functions $(x_1-\xi x_0)/x_0$ and $(x_3-\rho x_2)/x_0$.
Write 
\[
f(t)=\prod_{i=0}^{d-1}(t-\xi_i), \quad\quad g(t)=\prod_{j=0}^{d-1}(t-\rho_j),
\]
where all $\xi_j,\rho_j\in\bar k$.  
Let $L_{ij}$ be the projective line given by $x_1-\xi_i x_0=x_3-\rho_j x_2=0$.
Let $D$ be the subgroup of $\Div_{Z_{\bar k}}(X_{\bar k})$
generated by the lines $L_{ij}$ with $i\neq 0$, $j\neq 0$. It is easy to see that 
$\Div_{Z_{\bar k}}(X_{\bar k})=\mathrm{Im}(p_2)+(D\oplus\Z h)$.
Hence we get a surjective map $\phi\colon D\oplus\Z h\to\Lambda$.
It is enough to prove that $\phi$ is an isomorphism. To prove the injectivity of $\phi$
suppose that $a[h]+\sum_{i,j=1}^{d-1}a_{ij}[L_{ij}]=0$ for some $a_{ij}\in\Z$.

The intersection pairing $\langle [L_{ij}], [L_{rs}]\rangle$
equals $-(d-2)$ if $i=r$ and $j=s$, $0$ if $i\not=r$ and $j\not=s$, and $1$ in all other cases.
Pairing our sum with $L_{0,0}$ we get $a=0$. Next, pairing it with each $L_{0,s}$
and $L_{r,0}$ we see that $\sum_{i=1}^{d-1}a_{is}=\sum_{j=1}^{d-1}a_{rj}=0$ for all $s$ and $r$. Finally, pairing with each $L_{rs}$, where $1\leq r,s\leq d-1$, we obtain $a_{rs}=0$. This
proves that $\phi$ is injective, hence an isomorphism.
\eprf

\bthe \label{Lambda}
Let $L_{f,g}$ be the compositum of $k(\mu_d)$ and 
the field of definition of $\Pic(Y_{\bar k})$, where $Y$ is given by
$(\ref{X})$ with $a=1$.
If $[L_{f,g}(\sqrt[d]{a}):L_{f,g}]=d$, then
the natural inclusion of $\Ga_k$-modules
$\Lambda\to\Pic(X_{\bar k})$ gives rise to an isomorphism 
\[
\H^1(k,\Lambda)\tilde\lra\H^1(k,\Pic(X_{\bar k})).
\]
\ethe

The assumption
$[L_{f,g}(\sqrt[d]{a}):L_{f,g}]=d$ holds when $a$ is general enough for
given $f$ and $g$. By Proposition \ref{Lfg} this condition follows from 
$[K_{f,g}(\sqrt[d]{a}):K_{f,g}]=d$, where $K_{f,g}$ 
is the compositum of $k(\mu_d)$, the splitting fields of $f$ and $g$, and the field
of definition of $\Hom(J(C_{f,\bar k}),J(C_{g,\bar k}))$, where $J(C_f)$ 
(respectively, $J(C_g)$) is the Jacobian
of the smooth plane curve given by $y^d=f(t)$ (respectively, by $y^d=g(t)$).

\medskip

\bprf[Proof of Theorem \ref{Lambda}] The proof consists of several steps. 

{\em Step} 1. Since $L_{f,g}$ contains $\Q(\mu_d)$, the extensions 
$L_{f,g}(\sqrt[d]{a})/k$ and $L_{f,g}(\sqrt[d]{a})/L_{f,g}$ are Galois. 
Write $G=\Gal(L_{f,g}(\sqrt[d]{a})/k)$ and 
$H=\Gal(L_{f,g}(\sqrt[d]{a})/L_{f,g})$. By assumption, $H\simeq\Z/d$.
The Galois group $\Gal(\bar k/L_{f,g}(\sqrt[d]{a}))$ acts trivially on 
$\Pic(X_{\bar k})$, thus $\Ga_k$ acts on $\Pic(X_{\bar k})$ via
its quotient $G$. Since $\Pic(X_{\bar k})$ is a free abelian group of finite rank 
\cite[Théorème XI.1.8]{SGA},
the inflation map $\H^1(G,\Pic(X_{\bar k}))\to \H^1(k,\Pic(X_{\bar k}))$ is an isomorphism.
To finish the proof it is enough to prove the following statements:

(a$_X$) $\H^1(H,\Pic(X_{\bar k}))=0$;
	
(b$_X$) $\Pic(X_{\bar k})^H=\Lambda$.

\medskip

{\em Step} 2. Let $Y$ be the surface given by (\ref{X}) with $a=1$.
We note that $H$ acts on $Y_{L_{f,g}}$ by $L_{f,g}$-automorphisms which multiply the coordinates
$x_2$ and $x_3$ by the same root of unity of degree $d$, and leave $x_0$ and $x_1$
invariant.

Absorbing $\sqrt[d]{a}$
into $x_2$ and $x_3$ defines an isomorphism of base changes of
$X$ and $Y$ to $L_{f,g}(\sqrt[d]{a})$. This isomorphism translates the action of $H$
on $\Pic(X_{\bar k})$ induced by the action of $\Ga_k$ on $\bar k$
into the action of $H$
on $\Pic(Y_{\bar k})$ induced by the action of $H$ by automorphisms of $Y_{L_{f,g}}$.
Thus (a$_X$) and (b$_X$) are respectively equivalent to

(a$_Y$) $\H^1(H,\Pic(Y_{\bar k}))=0$;
	
(b$_Y$) $\Pic(Y_{\bar k})^H=\Lambda$.

\medskip

{\em Step} 3. The surface $Y_{L_{f,g}}$ together with the action of $H$ on $Y_{L_{f,g}}$ by automorphisms
can be defined over a subfield $L_0\subset L_{f,g}$ that is finitely generated over $\Q$.
Choose an embedding of $L_0$ into $\C$, and call $Y_\C$ the resulting complex surface.
This gives an isomorphism of $H$-modules $\Pic(Y_{\bar k})\cong\Pic(Y_\C)$.
Since $\H^1(Y_\C,\mathcal O_{Y_\C})=0$,
the exponential sequence gives an exact sequence of $H$-modules
\begin{equation}
0\to \Pic(Y_\C)\to\H^2(Y_\C,\Z)\to T\to 0.
\label{exp}
\end{equation}
By the Lefschetz theorem on $(1,1)$-classes,
$\Pic(Y_\C)$ is a saturated sublattice of the free abelian group $\H^2(Y_\C,\Z)$,
so that the abelian group $T$ is also free of finite rank.
The associated exact sequence
of cohomology groups of $H$, in view of semisimplicity of finite-dimensional  
representations of $H$ over a field of characteristic zero, shows that (a$_Y$) and (b$_Y$) follow from 

(a$_\C$) $\H^1(H,\H^2(Y_\C,\Z))=0$;

(b$_\C$) $\H^2(Y_\C,\Z)^H=\Lambda$.

\noindent These properties are purely topological, so they are invariant under
a continuous real deformation of complex polynomials $f(t)$ and $g(t)$ without multiple roots.
Therefore, it is enough to establish these properties for any chosen $f(t)$ and $g(t)$.
We shall check that (a$_\C$) and (b$_\C$) hold for 
the complex Fermat surface $F\subset\P^3_\C$, defined in (\ref{Fermat}),
that is, in the case when $f(t)=g(t)=t^d+1$. 

\medskip

{\em Step} 4. We claim that $\H^1(H,\H^2(F,\Z))=0$.
Fix a primitive $2d$-th root of unity $\zeta$ and write $\e=\zeta^2$.
We define $G\simeq(\Z/d)^3$ as the group of automorphisms of $F$ 
with generators $u_1,u_2,u_3$, where each $u_i$ multiplies $x_i$ by $\e$
and leaves the other variables invariant. 
Thus $H$ is the subgroup of $G$ generated by $u_2u_3$.

It is enough to show that $\H^1(H,P)=0$, where 
$P$ is the orthogonal complement to the hyperplane section class
$[h]$ in $\H^2(F,\Z)$. This is immediate from the
short exact sequence of $G$-modules
\[
0\to P\to \H^2(F,\Z)\to \Z\to 0.
\]
Let $\phi(x)\in\Z[x]$ be the polynomial $\sum_{i=0}^{d-1} x^i$. By \cite[Lemma 2.1]{GS22}
(based on work of Pham) there is an isomorphism of $G$-modules
\[
P\cong\Z[u_1,u_2,u_3]/I, \ \text{where} \ I=(\phi(u_1), \phi(u_2), \phi(u_3), \phi(u_1u_2u_3)).
\]
The group $\H^1(H,P)$ is the quotient of the subgroup of $P$ annihilated by $\phi(u_2u_3)$,
by $(1-u_2u_3)P$. Multiplication by $\phi(u_2u_3)$ on $P/(1-u_2u_3)P$ coincides with 
multiplication by $d$. We note that 
\[
P/(1-u_2u_3)P\cong\Z[u_1,u_2]/(\phi(u_1),\phi(u_2))
\]
is a torsion-free abelian group, so we conclude that $\H^1(H,P)=0$. 

\medskip

{\em Step} 5. We claim that $\H^2(F_\C,\Z)^H=\Lambda$. 
The group $H$ acts trivially on $\Lambda$, hence $\Lambda\subset\H^2(F_\C,\Z)^H$.

For a character $\chi\in \widehat G=\Hom(G,\Q/\Z)$, let $V_\chi\subset P\otimes\C$ be the eigenspace on which $G$ acts via $\chi$, i.e.\ such that $g(v)=\chi(g)v$ for any $g\in G$ and $v\in V_\chi$. By work of Shioda \cite{Sh79,Sh81}, see also \cite[Lemma 2.6]{GS22}, $V_\chi$ is non-zero, or equivalently $1$-dimensional, if and only if $\chi$ can be given by 
$\chi(u_1)=l_1/d$, $\chi(u_2)=l_2/d$, $\chi(u_3)=l_3/d$,
where $1\leq l_1,l_2,l_3\leq d-1$ and $l_1+l_2+l_3$ is not divisible by $d$.
Thus $P^H\otimes\C$ is the direct sum of the eigenspaces $V_\chi$ such that $l_2+l_3=d$.

The $d^2$ lines in $\Lambda$ can be written as follows:
\[L_{ij}:\quad \zeta x_1-\e^i x_0=x_3-\e^j\zeta x_2=0, \quad\quad i,j=0,\dots,d-1.\]
Following Shioda \cite[p.~732]{Sh81},
for $\chi=(l_1,d-l_3,l_3)\in \widehat G$ we define
\[\omega\chi=\sum_{i,j=0}^{d-1}\e^{il_1+jl_3}[L_{ij}]\in \Pic(F)\otimes\C.\]
Then $g\omega\chi=\chi(g)\omega\chi$ and 
$\langle \omega\chi,\ov{\omega\chi}\rangle=-d^3$,
hence $\omega\chi\neq 0$ spans $V_\chi$. 
Thus $V_\chi\subset\Lambda\otimes\C$, so that
$P^H\otimes\C\subset\Lambda\otimes\C$. Since $[h]\in\Lambda$,
we deduce that $\H^2(F_\C,\C)^H=\Lambda\otimes\C$. 

It follows that
$\Lambda$ is a subgroup of $\H^2(F_\C,\Z)^H$ of finite index.
In particular, $\Lambda$ has finite index in $\Pic(F)^H$.
By a result of Degtyarev \cite[Thm.~1.4]{Deg15}, the abelian group 
$\Pic(F)^H/\Lambda$ is torsion-free, thus $\Lambda=\Pic(F)^H$.
Since the abelian group $T$ in (\ref{exp}) is torsion-free, we must have
$\Pic(F)^H=\H^2(F_\C,\Z)^H$, which proves our claim. 
This finishes the proof of the theorem.
\eprf

\bco \label{2.4}
Assume that $[L_{f,g}(\sqrt[d]{a}):L_{f,g}]=d$. Then
the isomorphism of Theorem \ref{Lambda} identifies  the differential $d_{1,1}$, up to sign,
with the composition of the connecting map
 $\partial\colon\H^1(k,\Lambda)\to\H^3(k,\Z)$ attached to the bottom sequence of $(\ref{duo})$,
and $a_*\colon\H^3(k,\Z)\to \H^3(k,\bar k^\times)$ induced by the map sending $1\in \Z$ to
$a\in \bar k^\times$.
\eco
\bprf This follows from the
well known fact that, up to sign, $d_{1,1}$ is the connecting map attached to the top 
sequence of (\ref{duo}), see \cite[Prop.~5.4.5]{CTS21}.
\eprf

\medskip

It would be interesting to see if the deformation technique used to prove Theorem \ref{Lambda} can preserve sufficient arithmetic information in other settings where the Galois group acts by geometric automorphisms on a family of varieties (cf.\ Theorem \ref{friday}).

\section{Differential $d_{1,1}$} \label{diff}

Let $k_f$ and $k_g$ be the splitting fields of $f$ and $g$, respectively. 
Let $G_f=\Gal(k_f/k)$, $G_g=\Gal(k_g/k)$, and $G=\Gal(k_fk_g/k)$.
We denote the $G_f$-module $\Z[V_f]/\Z\sigma_f$ by $M_f$. Similarly, we have a $G_g$-module
$M_g=\Z[V_g]/\Z\sigma_g$. 

It is clear that the structure of a $\Ga_k$-module on $\Lambda$ is obtained by inflation
from its structure of a $G$-module.
Passing to the quotient by $\Z$ in the bottom sequence of (\ref{duo}) 
we obtain a commutative diagram of $G$-modules with exact rows
\begin{equation}
\begin{split}
\xymatrix{
&&&0\ar[d]&0\ar[d]&\\
&&&\Z\ar[d]\ar[r]^\id&\Z\ar[d]&\\
0\ar[r]& \Z\ar[r]^{p_1 \ \ \ \ \ \ \  \ }\ar[d]^{\id}& 
\Z[V_f]\oplus\Z[V_g]\ar[r]^{p_2 \ \ \ \ \ }\ar[d]^{\id}& 
(\Z[V_f]\otimes \Z[V_g])\oplus\Z h\ar[r]\ar[d]& \Lambda\ar[r]\ar[d]& 0\\
0\ar[r]& \Z\ar[r]^{p_1 \ \ \ \ \ \ \  \ }& 
\Z[V_f]\oplus\Z[V_g]\ar[r]& 
\Z[V_f]\otimes \Z[V_g]\ar[r]\ar[d]& M_f\otimes M_g\ar[r]\ar[d]& 0\\
&&&0&0&
}
\label{bruckner}
\end{split}
\end{equation}
The commutativity of (\ref{bruckner}) implies that 
$\partial\colon\H^1(G,\Lambda)\to \H^3(G,\Z)$ is the composition
\begin{equation}
\H^1(G,\Lambda)\hookrightarrow \H^1(G,M_f\otimes M_g)\stackrel{\partial_M}\lra \H^3(G,\Z),\label{new3}
\end{equation}
where $\partial_M$ is the connecting map attached to the bottom sequence of 
(\ref{bruckner}). In the next proposition, we 
calculate the group $\H^1(G,M_f\otimes M_g)$. Then we determine $\H^1(G,\Lambda)$
as a subgroup of this group. Finally, we compute the connecting map $\partial_M$.

For a group $G$ we write $G^*=\Hom(G,\Q/\Z)$.

\bpr \label{new1}
Assume that $f$ and $g$ are irreducible polynomials of the same degree
such that their splitting fields $k_f$ and $k_g$ are linearly disjoint over $k$. 
Let $S_f\subset G_f$ be the stabiliser of a $\bar k$-point of $V_f$, and let
$N_f\subset G_f$ be the normal closure of $S_f$ in $G_f$.
Define $S_g$ and $N_g$ similarly.
Then we have a natural isomorphism
\begin{equation}
\H^1(G, M_f\otimes M_g)\cong((G_f/N_f)^\mathrm{ab}\otimes (G_g/N_g)^\mathrm{ab})^*.\label{new5}
\end{equation}
If $\deg(f)=\deg(g)=d$, then 
$\H^1(G, M_f\otimes M_g)\cong\Hom((G_f/N_f)^\mathrm{ab}\otimes (G_g/N_g)^\mathrm{ab},\Z/d)$.
\epr
\bprf  The exact sequence of $G_f$-modules
$0\to\Z\to\Z[V_f]\to M_f\to 0 $
gives rise to the exact sequence
\[
0\to\Z\tilde\lra \Z[V_f]^{G_f}\to M_f^{G_f}\to \H^1(G_f,\Z)=0,
\]
where the second arrow is an isomorphism since $f$ is irreducible.
Thus $M_f^{G_f}=0$. Similarly, $M_g^{G_g}=0$.
We also see that $\H^1(G_f,M_f)$ is the kernel of the restriction map
$G_f^*\to S_f^*$. We deduce an isomorphism
\begin{equation}
\H^1(G_f,M_f)=((G_f/N_f)^\mathrm{ab})^*,\label{new2}
\end{equation}
and a similar isomorphism for $\H^1(G_g,M_g)$.

We have $G=G_f\times G_g$.
The K\"unneth formula for group cohomology in degree 1 gives the following (split) exact
sequence of abelian groups:
\[
0\to\bigoplus_{i=0,1}\big( \H^i(G_f,M_f)\otimes \H^{1-i}(G_g,M_g)\big)
\to\H^1(G, M_f\otimes M_g)
\]
\[
\to\bigoplus_{i=0,1,2}\Tor^\Z_1(\H^i(G_f,M_f),\H^{2-i}(G_g,M_g))\to 0.
\]
Since  $M_f^{G_f}=0$ and $M_g^{G_g}=0$, this gives an isomorphism
\[
\H^1(G, M_f\otimes M_g) \cong
\Tor^\Z_1(\H^1(G_f,M_f),\H^1(G_g,M_g)).
\]
An easy calculation shows that for finite abelian groups $A$ and $B$
there is a canonical isomorphism $\Tor^\Z_1(A^*,B^*)\cong (A\otimes B)^*$. Combining these results proves (\ref{new5}). 

If $\deg(f)=d$, then the index of $N_f$ in $G_f$ divides $d=[G_f:S_f]$. 
This implies the last statement of the proposition.
\eprf

\brem{
Let $T$ be the torus given by $N_{k(\xi)/k}(\mathbf{x})=N_{k(\rho)/k}(\mathbf{y})\neq 0$.
The pointed affine cone over our surface (\ref{X}) is a closed subset
of the smooth locus $S$ of the affine hypersurface $N_{k(\xi)/k}(\mathbf{x})=aN_{k(\rho)/k}(\mathbf{y})$,
which is a partial compactification of a $k$-torsor for $T$. The resulting
isomorphism of Galois modules $\Pic(S_{\bar k})\cong M_f\otimes M_g$
explains the similarity between Proposition \ref{new1} and \cite[Cor.~3.2]{Wei}.
}
\erem

\bpr \label{new4}
Assume that $f$ and $g$ are irreducible polynomials of degree $d$
such that their splitting fields $k_f$ and $k_g$ are linearly disjoint over $k$. 
The inclusion in $(\ref{new3})$ and the isomorphism $(\ref{new5})$ give rise to an exact sequence
\begin{equation}
0\to\H^1(G,\Lambda)\to
\Hom((G_f/N_f)^\mathrm{ab}\otimes (G_g/N_g)^\mathrm{ab},\Z/d)\stackrel{\chi}\lra
\Hom(G,\Z/d).
\label{new6}
\end{equation}
Here $\chi$ sends a bilinear function 
$\varphi\colon (G_f/N_f)^\mathrm{ab}\otimes (G_g/N_g)^\mathrm{ab}\to \Z/d$
to the homomorphism $\chi_\varphi\colon G\to\Z/d$ whose value on $(a,b)$, with
$a\in G_f$ and $b\in G_g$, is
\[
\chi_\varphi(a,b)=\sum_{y\in G_g/S_g}\widetilde\varphi(a,y)+\sum_{x\in G_f/S_f}\widetilde\varphi(x,b).
\]
In particular, $\H^1(G,\Lambda)\cong\Pi_{f,g}$.
\epr
\bprf In view of (\ref{new3}), $d$ annihilates $\H^1(G,\Lambda)$ since it annihilates
$\H^1(G, M_f\otimes M_g)$. Since $M_f^{G_f}=0$, the exact sequence
\[
0\to M_f\otimes M_g\stackrel{[d]}\lra M_f\otimes M_g\to (M_f\otimes M_g)/d\to 0
\]
gives rise to an isomorphism 
\begin{equation}
\label{may}
\big((M_f\otimes M_g)/d\big)^{G}\cong \H^1(G, M_f\otimes M_g).
\end{equation}
Likewise, by the right hand column of (\ref{bruckner}) we have $\Lambda^{G}=\Z[h]$, so
the exact sequence
\[
0\to\Lambda\stackrel{[d]}\lra\Lambda\to\Lambda/d\to 0
\]
gives rise to an isomorphism
$(\Lambda/d)^{G}/(\Z/d)[h]\cong \H^1(G, \Lambda)$.
Thus we need to determine the $G$-invariant elements of $(M_f\otimes M_g)/d$ which
lift to $(\Lambda/d)^G$.

The elements of $(M_f\otimes M_g)/d$ bijectively correspond to 
the functions $\varphi\colon G\to\Z/d$
that are constant on the left cosets by $S_f\times S_g$, considered modulo sums of functions $G_f/S_f\to\Z/d$ and $G_g/S_g\to\Z/d$. We can eliminate this indeterminacy by
normalising $\varphi$ so that $\varphi(e,y)=\varphi(x,e)=0$ for any $x\in G_f$ and $y\in G_g$.
The action of $h\in G_f$ sends the normalised $\varphi(x,y)$ to 
$\varphi(h^{-1}x,y)-\varphi(h^{-1},y)$. Thus a normalised function $\varphi$ is 
$G$-invariant if and only if it is bilinear. Equivalently, $\varphi$
is a bilinear function $(G_f/N_f)^\mathrm{ab}\otimes (G_g/N_g)^\mathrm{ab}\to\Z/d$.
(This gives another proof of the isomorphism (\ref{new5}).)

The elements of $\Lambda/d$ bijectively correspond to 
the functions $\psi\colon G\to\Z/d$
that are constant on the left cosets by $S_f\times S_g$, considered modulo sums
$\sigma(x)+\tau(y)$, where 
functions $\sigma\colon G_f/S_f\to\Z/d$ and $\tau\colon G_g/S_g\to\Z/d$ are
such that
\[
\sum_{x\in G_f/S_f}\si(x)+\sum_{y\in G_g/S_g}\tau(y)=0.
\]

Take an element of $((M_f\otimes M_g)/d)^G$ and represent it by
a function $\varphi\colon G\to\Z/d$ coming from a bilinear function 
$(G_f/N_f)^\mathrm{ab}\otimes (G_g/N_g)^\mathrm{ab}\to\Z/d$. 
It lifts to an element of $\Lambda/d$ uniquely up to an
element of $(\Z/d)[h]$, so one lift is $G$-invariant if and only if every lift is.
Thus it is enough to consider the lift given by the same function $\varphi$. 
Let $a\in G_f$ and $b\in G_g$.
In $\Lambda/d$ we have
\[
(a,b)\varphi(x,y)=\varphi(a^{-1}x,b^{-1}y)=
\varphi(x,y)+\varphi(x,b^{-1})+\varphi(a^{-1},y)+\varphi(a^{-1},b^{-1}).
\]
Thus $\varphi$ represents an element of $(\Lambda/d)^G$ if and only if the character 
$\chi_\varphi\colon G\to\Z/d$ given by
\[
\chi_\varphi\big((a,b)\big):=
\sum_{y\in G_g/S_g}\varphi(a^{-1},y)+\sum_{x\in G_f/S_f}\varphi(x,b^{-1})
\]
is trivial.
\eprf

\brem \label{bc}
{
Suppose that $f$ and $g$ are cyclic polynomials whose splitting fields are linearly disjoint over $k$.
We have $G_f\simeq G_g\simeq\Z/d$, $N_f=S_f=\{e\}$, 
$N_g=S_g=\{e\}$. Then (\ref{new6}) becomes the exact sequence
\[
0\to\H^1((\Z/d)^2,\Lambda)\to\Z/d\to \Hom((\Z/d)^2,\Z/d).
\]
The generator $1\in\Z/d$ goes to the character $\chi\colon(\Z/d)^2\to\Z/d$ 
such that $\chi((1,0))=\chi((0,1))$ is the sum of all elements of $\Z/d$.
This is zero if $d$ is odd, and  $d/2\in\Z/d$ if $d$ is even.
Thus in this case we have $\H^1(k,\Lambda)\cong 2(\Z/d)$. }
\erem

The Künneth formula for group cohomology gives an exact sequence
\[
0\to\H^3(G_f,\Z)\oplus \H^3(G_g,\Z)\to\H^3(G,\Z)\to
\Tor_1^\Z(\H^2(G_f,\Z),\H^2(G_g,\Z))\to 0,
\]
see \cite[Exercise 6.1.8]{Weibel}. The restriction
map $\H^3(G,\Z)\to \H^3(G_f,\Z)\oplus \H^3(G_g,\Z)$ splits this sequence, so we get
an isomorphism $\H^3(G,\Z)_\mathrm{prim}\cong(G_f^\mathrm{ab}\otimes G_g^\mathrm{ab})^*$.

\bpr \label{zoe}
The connecting map $\partial_M\colon\H^1(G,M_f\otimes M_g)\to \H^3(G,\Z)$ 
is the composition
\[
((G_f/N_f)^\mathrm{ab}\otimes (G_g/N_g)^\mathrm{ab})^*\hookrightarrow
(G_f^\mathrm{ab}\otimes G_g^\mathrm{ab})^*\tilde\lra \H^3(G,\Z)_\mathrm{prim}
\hookrightarrow \H^3(G,\Z).
\]
\epr
\bprf A choice of a $\bar k$-point in $V_f$ defines a surjective map of sets
$G_f\to V_f$ compatible with the (left) action of $G_f$.
The induced map of $G_f$-modules $\Z[V_f]\to\Z[G_f]$ and a similar map for $g$
give rise to a commutative diagram of $G$-modules
\begin{equation}\label{sixterms}\begin{split}
 \xymatrix{
0\ar[r]& \Z\ar[r]^{p_1 \ \ \ \ \ \ \  \ }\ar[d]^\id& \Z[V_f]\oplus\Z[V_g]\ar[r]\ar[d]& 
\Z[V_f]\otimes \Z[V_g]\ar[r]\ar[d]& M_f\otimes M_g\ar[r]\ar[d]^\kappa& 0\\
0\ar[r]& \Z\ar[r]& \Z[G_f]\oplus\Z[G_g]\ar[r]& 
\Z[G_f]\otimes \Z[G_g]\ar[r]&\widetilde M_f\otimes \widetilde M_g\ar[r]& 0}
\end{split}
\end{equation}
where $\widetilde M_f=\Z[G_f]/\Z$ and $\widetilde M_g=\Z[G_f]/\Z$.
We deduce that $\partial_M$ is the composition 
\[
\H^1(G,M_f\otimes M_g)\stackrel{\kappa_*}\lra \H^1(G,\widetilde M_f\otimes \widetilde M_g)
\stackrel{\partial_{\widetilde M}}\lra \H^3(G,\Z),
\]
where $\partial_{\widetilde M}$ is the connecting map of the
bottom exact sequence.  Proposition \ref{new1} gives an isomorphism
$\H^1(G,\widetilde M_f\otimes \widetilde M_g)\cong(G_f^\mathrm{ab}\otimes G_g^\mathrm{ab})^*$.
The natural map from
$G_f^\mathrm{ab}\otimes G_g^\mathrm{ab}$ to $(G_f/N_f)^\mathrm{ab}\otimes (G_g/N_g)^\mathrm{ab}$
is surjective, hence $\kappa_*$ is injective.

Let us write $M=(\Z[G_f]\oplus\Z[G_g])/\Z$. Since
$\Z[G_f]\otimes\Z[G_g]\cong\Z[G]$ is a free $\Z[G]$-module, we have a canonical
isomorphism $\H^1(G,\widetilde M_f\otimes \widetilde M_g)\cong\H^2(G,M)$.
The map \[\H^2(G,\Z)\to \H^2(G,\Z[G_f]\oplus\Z[G_g])\] is the natural isomorphism
$G^*\tilde\lra G_g^*\oplus G_f^*$, thus the connecting map of the exact sequence of $G$-modules
\[
0\to\Z\to\Z[G_f]\oplus\Z[G_g]\to M\to 0
\]
sends $\H^2(G,M)$ isomorphically onto $\H^3(G,\Z)_\mathrm{prim}$. 

It remains to check that the resulting isomorphism 
$(G_f^\mathrm{ab}\otimes G_g^\mathrm{ab})^*\tilde\lra \H^3(G,\Z)_\mathrm{prim}$
agrees with the isomorphism provided by the Künneth formula.

The group $\H^i(G_f,\Z)$, $i\geq 0$, is the $i$-th cohomology group of the complex
$C_f^\bullet$ of abelian groups where $C_f^i$ is the group of functions $(G_f)^i\to \Z$.
From the distinguished triangle $\Z\to\Z[G_f]\to \widetilde M_f\to\Z[1]$ it follows that $\mathrm{R}\Gamma(\widetilde M_f)\to C_f^\bullet[1]$ is a quasi-isomorphism in degrees $\geq 0$,
so that $\H^i(G_f,\widetilde M_f)$ is the $i$-th cohomology group of 
$(\tau_{\geq 1}C_f^\bullet)[1]$. 
The differential $C_f^0\to C_f^1$ sends $n\in\Z$ to $(g-1)n=0$
since $\Z$ is a trivial $G$-module. Thus $C_f^\bullet$ is the direct sum
of $\Z$ in degree 0 and the truncated complex $\tau_{\geq 1}C_f^\bullet$.
This implies that for $i\geq 2$ we have 
\[
\H^i(G,\Z)=\H^i(C_f^\bullet\otimes C_g^\bullet)=\H^i(C_f^\bullet)\oplus
\H^i(C^\bullet_g)\oplus\H^{i-2}\big((\tau_{\geq 1}C_f^\bullet)[1]\otimes 
(\tau_{\geq 1}C_g^\bullet)[1]\big),
\]
so the K\"unneth formula isomorphism is given by 
\[
\H^1\big((\tau_{\geq 1}C_f^\bullet)[1]\otimes 
(\tau_{\geq 1}C_g^\bullet)[1]\big)\to \H^3(C_f^\bullet\otimes C_g^\bullet).
\]
On the other hand, we can notice that the bottom 2-extension of the previous diagram (\ref{sixterms})
is equivalent to the Yoneda product of $0\to\Z\to\Z[G_f]\to\widetilde M_f\to 0$ and 
\[
0\to \widetilde M_f\to \widetilde M_f\otimes\Z[G_g]\to 
\widetilde M_f\otimes \widetilde M_g\to 0.
\]
Thus the resulting connecting map is the composition 
\[
\H^1\big((\tau_{\geq 1}C_f^\bullet)[1]\otimes (\tau_{\geq 1}C_g^\bullet)[1]\big)\to 
\H^2\big((\tau_{\geq 1}C_f^\bullet)[1]\otimes C_g^\bullet\big)
\to \H^3(C_f^\bullet\otimes C_g^\bullet),
\]
and hence coincides with the K\"unneth formula isomorphism. \eprf

\brem \label{bc1}
{This continues Remark \ref{bc} and uses the same assumptions.

(1) The cup-product map
\begin{equation}
\H^1(G_f,\Z/d)\otimes\H^1(G_g,\Z/d)\to \H^2(G,\Z/d)_\mathrm{prim}\label{cupp}
\end{equation}
is surjective. Indeed,
it is clear by functoriality that the image of the cup-product is contained in
$\H^2(G,\Z/d)_\mathrm{prim}$. Arguing exactly as in the topological
case (see diagram (5.37) in \cite[\S 5.7.2]{CTS21}) we identify the map (\ref{cupp}) with the map
\[
\Hom(G_f^\mathrm{ab},\Z/d)\otimes\Hom(G_g^\mathrm{ab},\Z/d)\to 
\Hom(G_f^\mathrm{ab}\otimes G_g^\mathrm{ab},\Z/d),
\]
given by multiplication in $\Z/d$. In our case $G_g^\mathrm{ab}\simeq\Z/d$, so
this is an isomorphism. 

(2)
The composition of the embedding $\H^2(G,\Z/d)_\mathrm{prim}\hookrightarrow \H^2(G,\Z/d)$ with
the connecting map $\delta\colon\H^2(G,\Z/d)\to \H^3(G,\Z)$ is an isomorphism.
Indeed, we have $\H^3(\Z/d,\Z)=0$, hence
$\H^3(G,\Z)_\mathrm{prim}=\H^3(G,\Z)\simeq\Z/d$, 
where the last isomorphism follows immediately from
the K\"unneth formula for group cohomology \cite[Exercise 6.1.8]{Weibel}.
Since $d$ annihilates $\H^2(G,\Z)\simeq(\Z/d)^2$, we have an exact sequence
\[
0\to \H^2(G,\Z)\to \H^2(G,\Z/d)\to \H^3(G,\Z)\to 0.
\]
Comparing it with similar sequences for $G_f$ and $G_g$ we see that the second arrow
can be identified with the canonical map $\H^2(G_f,\Z/d)\oplus\H^2(G_g,\Z/d)\to \H^2(G,\Z/d)$.
We conclude that $\delta\colon\H^2(G,\Z/d)_\mathrm{prim}\to \H^3(G,\Z)$ is an isomorphism.

(3)
Combining the results of (1) and (2), we see that the image of
the map $\H^3(G,\Z)\to \H^3(k,\Z)$ is generated by $\inf\delta(\chi_f\cup\chi_g)$,
where $\chi_f$ is a generator of $\H^1(G_f,\Z/d)\simeq\Z/d$ and $\chi_g$ is a generator of 
$\H^1(G_g,\Z/d)\simeq\Z/d$.
By Proposition \ref{zoe} the group $\H^1(G,M_f\otimes M_g)\cong \Z/d$ maps isomorphically 
onto $\H^3(G,\Z)$. Thus the image of $\H^1(G,M_f\otimes M_g)$ in $\H^3(k,\Z)$
is generated by $\inf\delta(\chi_f\cup\chi_g)$.

(4) 
Suppose that $\Q(\mu_{d})\subset k$. 
Fixing a primitive $d$-th root of unity, we identify $\Z/d$ with $\mu_d$. 
Write $f(t)=t^d-b$, $g(t)=t^d-c$, where $b,c\in k^\times$. Recall that
for $a\in k^\times$ we denote by $[a]\in\H^1(k,\mu_d)$ the image of $a$ under the
connecting map $k^\times\to\H^1(k,\mu_d)$ of the Kummer exact sequence.
We have $[b]=\inf\chi_f$ and $[c]=\inf\chi_g$.

We claim that the image of 
\[
\H^1(G,M_f\otimes M_g)\to\H^3(k,\Z)\stackrel{a}\lra\H^3(k,\bar k^\times),
\]
where the map $\Z\to\bar k^\times$ sends $1$ to $a$,
is generated by the image of $[a]\cup[b]\cup[c]$ under the map  
$\H^3(k,\mu_d^{\otimes 3})\to \H^3(k,\bar k^\times)$.
Indeed, there is a commutative diagram
\begin{equation}\begin{split}
\xymatrix{\H^2(k,\Z/d)\ar[r]^\delta\ar[d]_{\cup[a]}&\H^3(k,\Z)\ar[d]^{a}\\
\H^3(k,\mu_d)\ar[r]&\H^3(k,\bar k^\times)}
\label{difficult}\end{split}
\end{equation}
This is proved exactly as in \cite[\S 1.4.4]{CTS21}.
This implies the claim.

The diagonal surface (\ref{eqX}) with coefficients $a_1,a_2,a_3$ is isomorphic to
the surface (\ref{X}) with $a=-a_2$, $f(t)=t^d-b$ with $b=-a_1$, and $g(t)=t^d-c$
with $c=-a_3/a_2$. We have
\[ [a]\cup[b]\cup[c]=[-a_2]\cup[-a_1]\cup[-a_3/a_2]=-[-a_1]\cup[-a_2]\cup[-a_3],\]
so the connecting map $\partial$ sends
$\H^1(k,\Lambda)\simeq 2(\Z/d)$ to the subgroup of $\H^3(k,\bar k^\times)$ generated by
the image of $2[-a_1]\cup[-a_2]\cup[-a_3]$. If condition $(\star)$ holds, then
$d_{1,1}$ sends $\H^1(k,\Pic(X_{\bar k}))\simeq 2(\Z/d)$ to
the subgroup generated by the image of $2[-a_1]\cup[-a_2]\cup[-a_3]$. Furthermore,
if $(\star\star)$ also holds, then $d_{1,1}$ is injective.
}
\erem

\section{Proofs of the main results} \label{proofs}

We now prove the main technical result of this paper.

\bthe \label{f=g}
Let $f$ and $g$ be irreducible monic polynomials in $k[t]$ of degree $d$
whose splitting fields are linearly disjoint over $k$.
Let $X$ be the surface $(\ref{X})$, where $a\in k^\times$ is such that 
$[L_{f,g}(\sqrt[d]{a}):L_{f,g}]=d$. 
Then the following statements hold.

\textrm{(i)} We have $\H^1(k,\Pic(X_{\bar k}))\cong \Pi_{f,g}$.

\textrm{(ii)} The differential $d_{1,1}\colon \H^1(k,\Pic(X_{\bar k}))\to\H^3(k,\bar k^\times)$
can be identified with $\Delta_a$,
hence there is an isomorphism $\Br_1(X)/\Br_0(X)\cong \Ker(\Delta_a)$.
\ethe
\bprf
(i) This follows from Theorem \ref{Lambda} and Proposition \ref{new4}.

(ii) Corollary \ref{2.4}, formula (\ref{new3}) and Proposition \ref{zoe} imply that
$d_{1,1}$ can be identified with the composition
\[
\Pi_{f,g} \hookrightarrow\Hom((G_f/N_f)^\mathrm{ab}\otimes (G_g/N_g)^\mathrm{ab},\Z/d)\hookrightarrow
(G_f^\mathrm{ab}\otimes G_g^\mathrm{ab})^*
\]
\[
\tilde\lra\H^3(G,\Z)_\mathrm{prim}\hookrightarrow\H^3(G,\Z)\stackrel{\mathrm{inf}}\lra
\H^3(k,\Z)\stackrel{a}\lra\H^3(k,\bar k^\times),
\]
where the last arrow is induced by the map sending $1\in\Z$ to $a\in \bar k^\times$.
The resulting map coincides with $\Delta_a$ defined in the introduction, as
follows from the commutativity of diagram (\ref{difficult}). \eprf

\medskip

It is noteworthy this proof uses neither any knowledge of 
the group $\Pic(X_{\bar k})$ (like an integral basis or the Galois action), nor of the field $L_{f,g}$.

\bprf[Proof of Proposition \ref{Pi}]

A transitive subgroup $H\subset S_d$ is primitive
if the stabiliser of a point is a maximal subgroup of $H$.
Assume that $G_f\subset S_d$ is primitive.
If $\Pi_{f,g}\neq 0$, then $N_f\neq G_f$, but since
$S_f$ is a maximal subgroup of $G_f$, the inclusion $S_f\subset N_f$ is an equality, hence
$S_f$ is normal in $G_f$.
In this case $S_f$ stabilises every root of $f$, hence $S_f=\{e\}$ and $|G_f|=d$.
If $d$ is not prime, then the trivial subgroup cannot be maximal in $G_f$, a contradiction.
If $d$ is prime, we get $G_f\simeq \Z/d$. 
For a prime $d$, any transitive group is primitive. Thus a similar argument gives $G_g\simeq\Z/d$.
For any $d$, if $G_f\simeq G_g\simeq \Z/d$, then $\Pi_{f,g}\simeq\Z/d$ if $d$ is odd
and $\Pi_{f,g}\simeq\Z/(d/2)$ if $d$ is even.
\eprf

\bprf[Proof of Theorem \ref{general}]

By Remark \ref{Hilb} we can find a $k$-point $u$ in $\P^3_k$
such that $F(u)\neq 0$, $G(u)\neq 0$ and the class of $G(u)/F(u)$ in $K^\times/K^{\times d}$
has order $d$.
After a linear change of variables $x_0,x_1$ (respectively, $x_2,x_3$)
we can assume that $u=(0:1:0:1)$. We can write
$F(x_0,x_1)=F_0 f(x_1/x_0)x_0^d$ and $G(x_2,x_3)=G_0 g(x_3/x_2)x_2^d$, where
the polynomials $f(t)$ and $g(t)$ are monic and $F_0, G_0\in k^\times$. 
We have $[G(u)/F(u)]=[G_0/F_0]$. Thus $X$
is given by equation (\ref{X}) with $a=G_0/F_0$. 
By Remark \ref{defi} we have $L_{f,g}\subset K$, so the assumption that
the class of $a$ in $K^\times/K^{\times d}$ has order $d$ implies
$[L_{f,g}(\sqrt[d]{a}):L_{f,g}]=d$. It remains to apply Theorem \ref{f=g}.
\eprf

\bprf[Proof of Theorem \ref{main} (i)] 

We consider the diagonal surface $X$ over $k$ given by (\ref{eqX}).
This is a particular case of the surface (\ref{X}) with $a=-a_2$,
\[
f(t)=t^d+a_1=\prod_{i=0}^{d-1}(t-\e^i\xi), \quad
g(t)=t^d+a_3/a_2=\prod_{j=0}^{d-1}(t-\e^j\rho),
\]
where $\xi,\rho\in\bar k$ are such that $\xi^d+a_1=\rho^d+a_3/a_2=0$.
We have $L_{f,g}\subset k_L(\xi,\rho)$, where
$L$ is the field of definition of the geometric Picard group of the Fermat surface, and $k_L=kL$.
Condition $(\star)$ implies that we can use Theorem \ref{Lambda}. 
The action of $\Ga_k$ on $\Lambda$ factors through 
$\mathfrak G=\Gal(k(\mu_d,\xi,\rho)/k)$, hence we have
\[
\H^1(k,\Pic(X_{\bar k}))\cong\H^1(k,\Lambda)\cong\H^1(\mathfrak G,\Lambda).
\]
Condition $(\star)$ implies that
$G=\Gal(k(\mu_d,\xi,\rho)/k(\mu_d))\cong (\Z/d)^2$, which acts on $\xi$ and $\rho$ by
multiplication by $d$-th roots of unity. Thus $\mathfrak G=G\rtimes \fH$,
where 
\[
\fH=\Gal(k(\mu_d)/k)=\Gal(k(\mu_d,\xi,\rho)/k(\xi,\rho))\subset (\Z/d)^\times.
\]
Consider the semidirect product $G\rtimes(\Z/d)^\times$, where
$(\Z/d)^\times$ acts on $G=(\Z/d)^2$ by multiplication. 
The lines in $\Lambda$ are naturally indexed by the elements of $G$.
This equips $\Lambda$ with an action of $G\rtimes(\Z/d)^\times$ 
which extends the action of $\mathfrak G$.

The right hand column of diagram (\ref{bruckner}) is an exact sequence
of $G\rtimes(\Z/d)^\times$-modules, hence also of $\mathfrak G$-modules
\begin{equation}
0\to\Z\to\Lambda\to M_f\otimes M_g\to 0.\label{monday}
\end{equation}
Let $\Lambda_0\subset \Lambda$ be the subgroup generated by the lines 
indexed by the elements of $G\cong (\Z/d)^2$ with both coordinates non-zero.
(In other terms, these are the lines $x_1-\e^i\xi x_0=x_3-\e^j\rho x_2=0$ for $i\neq 0$, $j\neq 0$.)
As a $(\Z/d)^\times$-module, $\Lambda$ is isomorphic to the direct sum $\Z\oplus\Lambda_0$,
so that (\ref{monday}) splits as an exact sequence of $(\Z/d)^\times$-modules.

From (\ref{monday}) we get an exact sequence
\begin{equation}
0\to \Z\to\Lambda^G\to 0\to\H^1(G,\Lambda)\to \H^1(G, M_f\otimes M_g)\to \H^1(G,\Z/d), \label{rrr}
\end{equation}
where we used that $(M_f\otimes M_g)^G=0$. Thus we have $\Lambda^G=\Z$.
Taking $\fH$-invariants we obtain an exact sequence
\begin{equation}
0\to\H^1(G,\Lambda)^\fH\to \H^1(G, M_f\otimes M_g)^\fH\to \H^1(G,\Z/d), \label{grrr}
\end{equation}
If we denote the two factors of $G$ by $G_1\cong\Z/d$ and $G_2\cong\Z/d$, then
by Proposition \ref{new1} applied over $k(\mu_d)$ the group $\H^1(G,M_f\otimes M_g)$ is canonically isomorphic to
$\Hom(G_1\otimes G_2,\Z/d)\cong\Z/d$. Thus
$n\in (\Z/d)^\times$ acts on this group as multiplication by $n^{-2}$.
We conclude that $\H^1(G,\Lambda)$, 
which by Remark \ref{bc}, again applied over $k(\mu_d)$, is the subgroup $2(\Z/d)\subset\Z/d$, 
inherits the same action of $(\Z/d)^\times$.

Since $d$ annihilates $\H^1(G,M_f\otimes M_g)$ and $\H^1(G,\Lambda)$, 
it annihilates $\H^1(\mathfrak G,M_f\otimes M_g)\cong \H^1(G,M_f\otimes M_g)^\fH$
and $\H^1(\mathfrak G, \Lambda)\subset \H^1(G,\Lambda)$.
Thus we have isomorphisms
\[
\big((M_f/d)\otimes (M_g/d)\big)^{\mathfrak G}\tilde\lra \H^1(\mathfrak G,M_f\otimes M_g)
\cong \H^1(G,M_f\otimes M_g)^\fH
\]
and $(\Lambda/d)^{\mathfrak G}/(\Z/d)[h]\cong \H^1(\mathfrak G, \Lambda)$.

Finally, the reduction of (\ref{monday}) modulo $d$ gives an exact sequence
\begin{equation}
0\to\Z/d\to(\Lambda/d)^G\to \big((M_f/d)\otimes (M_g/d)\big)^G\to\H^1(G,\Z/d).\label{22}
\end{equation}
Since (\ref{monday}) splits as an exact sequence of $\fH$-modules, if $x\in\Lambda/d$
goes to an $\fH$-invariant element of $(M_f/d)\otimes (M_g/d)$, then $x\in(\Lambda/d)^\fH$.
Thus (\ref{22}) gives rise to an exact sequence
\begin{equation}
0\to\Z/d\to(\Lambda/d)^{\mathfrak G}\to \big((M_f/d)\otimes (M_g/d)\big)^{\mathfrak G}
\to\H^1(G,\Z/d).
\end{equation}
The cokernel of the second arrow, which is $\H^1(\mathfrak G, \Lambda)$,
 is isomorphic to the kernel of the last arrow, which is $\H^1(G,\Lambda)^{\mathfrak H}$, see
(\ref{grrr}). This proves that $\H^1(k,\Pic(X_{\bar k}))$ is the subgroup of
$2(\Z/d)$ on which the square of the cyclotomic character acts by $1$.
\eprf

\bprf[Proof of Theorem \ref{main} (ii)]
If $\H^1(k,\Pic(X_{\bar k}))=0$, the statement is clear from the spectral sequence (\ref{ss}).

The statement concerning the case $\Q(\mu_{d})\subset k$ follows from Remark \ref{bc1}.

Now assume that condition $(\star\star)$ holds. We have a commutative diagram
\[
\xymatrix{\H^1(k(\mu_{d}),\Lambda)\ar[r]&\H^3(k(\mu_{d}),\bar k^\times)\\
\H^1(k,\Lambda)\ar[r]\ar[u]&\H^3(k,\bar k^\times)\ar[u]}
\]
The left hand restriction map is injective: the 
$\Gal(\bar k/k(\mu_{d}))$-invariant subgroup of $\Lambda$ is generated by $[h]$,
but $\H^1(\Gal(k(\mu_{d})/k),\Z)=0$. By Remark \ref{bc1}
the top horizontal map sends $\H^1(k(\mu_{d}),\Lambda)\simeq 2(\Z/d)$ to
the subgroup generated by the image of $2[-a_1]\cup[-a_2]\cup[-a_3]$ in 
$\H^3(k(\mu_{d}),\bar k^\times)$, so by
condition $(\star\star)$ it is injective.
Thus the bottom horizontal map is also injective. By Theorem \ref{Lambda}
this implies that $d_{1,1}\colon\H^1(k,\Pic(X_{\bar k}))\to\H^3(k,\bar k^\times)$
is injective. The statement follows from the spectral sequence (\ref{ss}).

It remains to prove that condition $(\star\star)$ holds when $k=k_0(a_1,a_2,a_3)$ is
purely transcendental of degree 3 over $k_0$.
 For this we may assume that $k_0$ is algebraically closed. 

The discrete valuation of $k$ given by the variable $a_3$ gives rise to the residue
maps $\H^3(k,\bar k^\times)\to \H^2(k_0(a_1,a_2),\Q/\Z)$, see \cite[\S 6.3]{GSz},
and $\H^3(k,\mu_d)\to\H^2(k_0(a_1,a_2),\Z/d)$, see  \cite[\S 6.8]{GSz}.
By \cite[Prop.~6.8.9]{GSz} the following diagram {\em anticommutes}:
\[
\xymatrix{\H^3(k,\bar k^\times)\ar[r]& \H^2(k_0(a_1,a_2),\Q/\Z)\\
\H^3(k,\mu_d)\ar[r]\ar[u]&\H^2(k_0(a_1,a_2),\Z/d)\ar[u]}
\]
Since $k_0$ is algebraically closed of characteristic zero, we have an isomorphism 
of $\Ga_k$-modules $\Q/\Z\cong (\Q/\Z)(1)$, which induces an isomorphism
\[\H^2(k_0(a_1,a_2),\Q/\Z)\cong \H^2(k_0(a_1,a_2),\Q/\Z(1)).\] The cokernel of 
the natural inclusion $\Q/\Z(1)\to \bar k^\times$ is uniquely divisible,
hence is a $\Q$-vector space. This implies that the induced map
\[\H^2(k_0(a_1,a_2),\Q/\Z(1))\to \H^2(k_0(a_1,a_2),\bar k^\times)\] is an isomorphism.

Combining the previous diagram with a similar (anticommuting) 
diagram for the discrete valuation given by the variable $a_2$ we obtain
\[
\xymatrix{\H^3(k,\bar k^\times)\ar[r]& \H^2(k_0(a_1,a_2),\Q/\Z)\ar[r]^\cong&
\H^2(k_0(a_1,a_2),\bar k^\times)\ar[r]&\H^1(k_0(a_1),\Q/\Z)\\
\H^3(k,\mu_d)\ar[r]\ar[u]&\H^2(k_0(a_1,a_2),\Z/d)\ar[u]\ar[r]^\cong&
\H^2(k_0(a_1,a_2),\mu_d)\ar[u]\ar[r]&\H^1(k_0(a_1),\Z/d)\ar[u]}
\]
It is easy to see that the right hand vertical map is injective. Using \cite[Prop.~6.8.7]{GSz}
we obtain that, up to sign, the image of $[a_1]\cup[a_2]\cup[a_3]$ in $\H^2(k_0(a_1,a_2),\Z/d)$
is $[a_1]\cup[a_2]$, and that, again up to sign, $[a_1]\cup[a_2]$ goes to 
$[a_1]\in \H^1(k_0(a_1),\Z/d)$, which has order $d$. This clearly implies
that the order of the image of $[a_1]\cup[a_2]\cup[a_3]$ in $\H^3(k,\bar k^\times)$
is $d$, which finishes the proof.
\eprf

\bprf[Proof of Theorem \ref{main} (iii)]

Recall that $G\simeq(\Z/d)^3$ is a group of automorphisms of the Fermat surface $F$ 
with generators $u_1,u_2,u_3$, where $u_i$ multiplies $x_i$ by $\e$
and leaves the other variables invariant. 
Note that $(P\otimes\C)^G=0$ implies that $\H^2(F_\C,\Q)^G=\Q[h]$, 
where $[h]$ is the hyperplane class. The intersection index of $[h]$ with the class of a
line on $F_\C$ is $1$, hence $\H^2(F_\C,\Z)^G=\Z[h]$. In particular, we obtain
$\Pic(F_\C)^G=\Z[h]$. This implies $\Pic(X_{\bar k})^{\Ga_k}=\Z[h]$.
But $[h]$ is contained in the image of the natural map
$\Pic(X)\to\Pic(X_{\bar k})^{\Ga_k}$, thus this map is an isomorphism.
The spectral sequence (\ref{ss}) now gives the injectivity of $\Br(k)\to\Br(X)$.
\eprf

\brem \label{3.6}
{
One can apply the same technique to the cases not covered by Theorems \ref{f=g}
and \ref{main}. For example, suppose that $f$ and $g$ are irreducible polynomials
such that $k_f=k_g$ and $G_f=G_g\simeq\Z/d$. Thus $G=\Gal(k_fk_g/k)\simeq\Z/d$.
Breaking the bottom row of (\ref{duo}) into two short exact sequences
\begin{equation}
0\to \Z\to \Z[V_f]\oplus\Z[V_g]\to M\to 0,\quad
0\to M\to (\Z[V_f]\otimes \Z[V_g])\oplus\Z \to \Lambda\to 0 \label{dve}
\end{equation}
and using
that $\Z[V_f]$ and $\Z[V_g]$ are free $\Z[G]$-modules we obtain embeddings
\[
\H^1(G,\Lambda)\hookrightarrow\H^2(G,M)\hookrightarrow\H^3(G,\Z)\cong\H^1(\Z/d,\Z)=0.
\]
Assuming $[L_{f,g}(\sqrt[d]{a}):L_{f,g}]=d$ we obtain
$\H^1(k,\Pic(X_{\bar k}))=0$ by Theorem \ref{Lambda}.
}
\erem

\brem \label{3.7}
{
Another similar case is when $g$ is irreducible with $G_g\simeq\Z/d$, but $f$
is a product of linear factors. We have $G=\Gal(k_fk_g/k)\simeq\Z/d$.
In this case (\ref{duo}) simplifies to the following exact sequence of $G$-modules:
\begin{equation}
0\to \Z^{d-2}\to \Z[G]^{d-1}\to\Lambda\to 0.
\end{equation}

We deduce an isomorphism $\H^1(G,\Lambda)\simeq(\Z/d)^{d-2}$.
Assuming $[L_{f,g}(\sqrt[d]{a}):L_{f,g}]=d$ we obtain
$\H^1(k,\Pic(X_{\bar k}))\simeq(\Z/d)^{d-2}$ by Theorem \ref{Lambda}.
}
\erem

\section{Transcendental Brauer group} \label{3}

We start with a simple general lemma.

\ble \label{useful}
Let $X$ be a smooth, integral variety over a field $k$ of characteristic zero.
If for every prime $\ell$ and for every large integer $n$
the cycle class map induces an isomorphism 
$(\NS(X_{\bar k})/\ell^n)^{\Ga_k}\tilde\lra \H^2_\et(X_{\bar k},\mu_{\ell^n})^{\Ga_k}$,
then the map $\Br(X)\to\Br(X_{\bar k})$ is zero.
\ele
\bprf For every prime power $\ell^n$
the Kummer exact sequences for $X$ and $X_{\bar k}$ give rise to
the following commutative diagram  with exact rows
\[
\xymatrix{0\ar[r]&(\NS(X_{\bar k})/\ell^n)^{\Ga_k}\ar[r]& 
\H^2_\et(X_{\bar k},\mu_{\ell^n})^{\Ga_k}\ar[r]&
\Br(X_{\bar k})[\ell^n]^{\Ga_k}&\\
&&\H^2_\et(X,\mu_{\ell^n})\ar[r]\ar[u]&\Br(X)[\ell^n]\ar[u]\ar[r]&0}
\]
Since $X$ is smooth, $\Br(X)$ is a torsion group \cite[Lemma 3.5.3]{CTS21}. Thus it is enough to
prove that the right hand vertical map is zero for all $\ell$ and all large
enough integers $n$. By the commutativity of the diagram, this holds if
the second map in the top exact sequence is an isomorphism.
\eprf

\bthe \label{friday}
Let $k_0$ be an algebraically closed field of characteristic zero, and let
$k$ be a field extension of $k_0$. 
Let $k\subset K\subset\bar k$ be a finite Galois extension of $k$ with Galois group 
$\g=\Gal(K/k)$. Let $c\colon\Ga_k\to\g$ be the $1$-cocycle for the trivial action of $\Ga_k$ on $\g$
given by the natural surjective homomorphism $\Ga_k\to\g$. 

Let $Y$ be a smooth, projective, integral variety over $k_0$
such that $\H^3_\et(Y,\Z_\ell)$ is torsion-free for all primes $\ell$. 
Suppose that $\g$ acts on $Y$ by automorphisms so that for every prime $\ell$ we have

\textrm{(a)} $\H^1(\g,\H^2_\et(Y,\Z_\ell)/\mathrm{tors})=0$;

\textrm{(b)} the cycle class map 
$c_1\colon(\NS(Y)\otimes\Q_\ell)^\g\to\H^2_\et(Y,\Q_\ell)^\g$ is surjective.

\noindent Let $X$ be the $(K/k)$-twist of $Y_k$ by the $1$-cocycle $c\colon\Ga_k\to\g$.
Then $\Br(X)\to\Br(X_{\bar k})$ is the zero map.
\ethe
\bprf Since $\Br(X)$ is a torsion group, it is enough to prove the triviality of
the map of $\ell$-primary torsion subgroups $\Br(X)\{\ell\}\to\Br(X_{\bar k})\{\ell\}$ for all $\ell$.

From the definition of $X$ we see that there is an
isomorphism $X_{\bar k}\cong Y_{\bar k}$ under which
the action of $g\in\Ga_k$ on $X_{\bar k}$ is the same as the action of $g$
on $Y_{\bar k}$ followed by the action of $c(g)\in \g$. The same holds for the induced
action on cohomology: for any prime $\ell$ and any $n\geq 1$, $i\geq 0$ we have 
isomorphisms
$\H^i_\et(X_{\bar k},\mu_{\ell^n})\cong\H^i_\et(Y_{\bar k},\mu_{\ell^n})$
translating the natural action of $g$ on $\H^i_\et(X_{\bar k},\mu_{\ell^n})$ into the natural
action of $g$ on $\H^i_\et(Y_{\bar k},\mu_{\ell^n})$ followed by the action of $c(g)$. 

Since $k_0\subset \bar k$ is an extension of algebraically closed fields, and
$Y_{\bar k}$ is obtained from $Y$ by base change, we have a canonical
isomorphism $\H^i_\et(Y,\mu_{\ell^n})\to \H^i_\et(Y_{\bar k},\mu_{\ell^n})$, thus
$\H^i_\et(Y_{\bar k},\mu_{\ell^n})$ is a trivial $\Ga_k$-module. It follows that 
the $\Ga_k$-module structure on $\H^i_\et(X_{\bar k},\mu_{\ell^n})$ 
is obtained from the $\g$-module $\H^i_\et(Y,\mu_{\ell^n})$ by inflation.
Using the Kummer sequence we identify the inclusion of
$\Ga_k$-modules $\NS(X_{\bar k})/\ell^n\hookrightarrow \H^2_\et(X_{\bar k},\mu_{\ell^n})$ 
with the inclusion of $\g$-modules
$\NS(Y)/\ell^n\hookrightarrow \H^2_\et(Y,\mu_{\ell^n})$.
In view of Lemma \ref{useful}, it is enough to prove that the induced map
 of $\g$-invariant subgroups is an isomorphism.

Let us write $\NS(Y)_\ell:=\NS(Y)\otimes\Z_\ell$.
We can choose $n$ so that the torsion subgroup of $\H^2_\et(Y,\Z_\ell)$
is annihilated by $\ell^n$. Then the cycle class map and multiplication by $\ell^n$ 
give rise to a commutative diagram with exact rows
\[
\xymatrix{0\ar[r]& \H^2_\et(Y,\Z_\ell)/\mathrm{tors}\ar[r]& \H^2_\et(Y,\Z_\ell)\ar[r]&
\H^2_\et(Y,\mu_{\ell^n})\ar[r]&0\\
0\ar[r]&\NS(Y)_\ell/\mathrm{tors}\ar[r]\ar[u]&\NS(Y)_\ell\ar[r]\ar[u]&\NS(Y)/\ell^n\ar[r]\ar[u]&0
}
\]
where we used that $\H^3_\et(Y,\Z_\ell)_\mathrm{tors}=0$.
By condition (a) we obtain the following commutative diagram with exact rows
\[
\xymatrix{
0\ar[r]& \big(\H^2_\et(Y,\Z_\ell)/\mathrm{tors}\big)^\g\ar[r]&\H^2_\et(Y,\Z_\ell)^\g\ar[r]&
\H^2_\et(Y,\mu_{\ell^n})^\g\ar[r]& 0\\
0\ar[r]&\big(\NS(Y)_\ell/\mathrm{tors}\big)^\g\ar[r]\ar[u]&\big(\NS(Y)_\ell\big)^\g\ar[r]\ar[u]&
\big(\NS(Y)/\ell^n\big)^\g\ar[u]&
}
\]
The desired surjectivity of $(\NS(Y)/\ell^n)^\g\to \H^2_\et(Y,\mu_{\ell^n})^\g$
follows from the surjectivity of $(\NS(Y)_\ell)^\g \to \H^2_\et(Y,\Z_\ell)^\g$, so let us
establish this.

The cokernel of $\NS(Y)_\ell \to \H^2_\et(Y,\Z_\ell)$
is torsion-free \cite[XI, Lemme 1.9 (i)]{SGA}, thus we have a commutative diagram of $\g$-modules with exact rows
\[
\xymatrix{
0\ar[r]&\H^2_\et(Y,\Z_\ell)_\mathrm{tors}\ar[r]&\H^2_\et(Y,\Z_\ell)\ar[r]&
\H^2_\et(Y,\Z_\ell)/\mathrm{tors}\ar[r]&0\\
0\ar[r]&\NS(Y)\{\ell\}\ar[r]\ar[u]^{\cong}&\NS(Y)_\ell\ar[r]\ar[u]&
\NS(Y)_\ell/\mathrm{tors}\ar[r]\ar[u]&0
}
\]
This gives rise to the following commutative diagram with exact rows
\[
\xymatrix@C=1em{0\ar[r]&\H^2_\et(Y,\Z_\ell)_{\mathrm{tors}}^\g\ar[r]&
\H^2_\et(Y,\Z_\ell)^\g\ar[r]&
\big(\H^2_\et(Y,\Z_\ell)/\mathrm{tors}\big)^\g\ar[r]&
\H^1(\g,\H^2_\et(Y,\Z_\ell)_\mathrm{tors})\\
0\ar[r]&\NS(Y)\{\ell\}^\g\ar[r]\ar[u]^{\cong}&
\big(\NS(Y)_\ell\big)^\g\ar[r]\ar[u]&
\big(\NS(Y)_\ell/\mathrm{tors}\big)^\g\ar[r]\ar[u]&
\H^1(\g, \NS(Y)\{\ell\})\ar[u]^{\cong}
}
\]
Condition (b) implies the surjectivity of $\big(\NS(Y)_\ell/\mathrm{tors}\big)^\g\to 
\big(\H^2_\et(Y,\Z_\ell)/\mathrm{tors}\big)^\g$. Now
the surjectivity of 
$(\NS(Y)_\ell)^\g \to \H^2_\et(Y,\Z_\ell)^\g$
follows from the 5-Lemma.
\eprf

\bprf[Proof of Theorem \ref{n}]

For any extension $k\subset k'\subset\bar k$ we have a commutative diagram
\begin{equation}
\begin{split}
\xymatrix{\Br(X_{k'})\ar[r]&\Br(X_{\bar k})\\
\Br(X)\ar[r]\ar[u]&\Br(X_{\bar k})\ar[u]_{\mathrm{id}}}\label{factor}
\end{split}
\end{equation}
Thus, in order to prove that the map $\Br(X)\to\Br(X_{\bar k})$ is zero,
we can replace $k_0$ by its algebraic closure in $\bar k$, and so
assume without loss of generality that $k_0$ is algebraically closed.

Let $Y$ be the surface (\ref{X}) with $a=1$ over $k_0$. 
By \cite[Théorème XI.1.8]{SGA} we have $\Pic(Y)_\mathrm{tors}=0$, hence $\H^2_\et(Y,\Z_\ell)_\mathrm{tors}=0$
for all $\ell$. Poincar\'e duality implies that $\H^3_\et(Y,\Z_\ell)_\mathrm{tors}=0$
for all $\ell$. 
This allows us to
apply Theorem \ref{friday} to $Y$, the extension $k=k_0(a)$ of $k_0$ of transcendence degree 1,
and the group $\g=H=\Gal(k(\sqrt[d]{a})/k)$ acting on $Y$ as in
the proof of Theorem \ref{Lambda}. 
Indeed, conditions (a) and (b) hold, as follows from
comparison theorems between classical and
$\ell$-adic \'etale cohomology and statements $(\mathrm{a}_\C)$ and $(\mathrm{b}_\C)$
in the proof of Theorem \ref{Lambda}. 

The last statement follows from Theorem \ref{f=g} (ii).
\eprf

\bprf[Proof of Theorem \ref{co1}]

We are in a particular case of the situation considered in Theorem \ref{n}, so we can apply
this theorem to deduce that the map $\Br(X)\to\Br(X_{\bar k})$ is zero.
It remains to apply Theorem \ref{main} (ii) and (iii).
\eprf

\subsubsection*{Further applications}

Our methods can be applied to other, `less generic' diagonal surfaces. Here we present several
examples, leaving it to the interested reader to explore other cases.

Let $k$ be a field containing $L$, the field of definition of
the geometric Picard group of the Fermat surface of degree $d$ over $\Q$. 
In particular, $\Q(\mu_{2d})\subset k$. 

\medskip

\noindent\textbf{Example 1} Let $t\in k^\times$ be such that the class of $t$ in 
$k^\times/k^{\times d}$ has order $d$.
Let $X\subset\P^3_k$ be the surface
\[
x_0^d+x_1^d+tx_2^d+tx_3^d=0.
\]
Theorem \ref{Lambda} implies that $\H^1(k,\Pic(X_{\bar k}))=0$. 
If there is a subfield $k_0\subset k$ such that $k=k_0(t)$ has transcendence degree 1 over $k_0$,
then $\Br(X)=\Br_1(X)$ by 
Theorem \ref{n}. Since $X$ has a $k$-point,
the map $\Br(k)\to\Br(X)$ is an isomorphism.

\medskip

\noindent\textbf{Example 2} Let $s,t\in k^\times$ be such that the subgroup of
$k^\times/k^{\times d}$ generated by the classes of $s$ and $t$ is isomorphic to $(\Z/d)^2$.
Let $X\subset\P^3_k$ be the surface
\[
x_0^d+sx_1^d+tx_2^d+s^rtx_3^d=0,
\]
where $r$ is an integer coprime to $d$.
Remark \ref{3.6} gives $\H^1(k,\Pic(X_{\bar k}))=0$. 
If there is a subfield $k_0\subset k$ such that $k=k_0(t)$ has transcendence degree 1 over $k_0$,
then $\Br(X)=\Br_1(X)$ by Theorem \ref{n}, thus $\Br(X)=\Br_0(X)$. 
Now assume that $r=1$ and there is a subfield $k_1\subset k$ such that $k=k_1(s,t)$
has transcendence degree 2 over $k_1$. In this case,
by Proposition \ref{suppl} below, the kernel of
the map $\Br(k)\to\Br(X)$ is generated by
the class of the cyclic algebra $(s,t)_d$, so is isomorphic to $\Z/d$.

\medskip

\noindent\textbf{Example 3} Let $s,t\in k^\times$ be such that the subgroup of
$k^\times/k^{\times d}$ generated by the classes of $s$ and $t$ is isomorphic to $(\Z/d)^2$. 
Let $X\subset\P^3_k$ be the surface
$$x_0^d+x_1^d+tx_2^d+stx_3^d=0.$$
Remark \ref{3.7} gives $\H^1(k,\Pic(X_{\bar k}))\simeq(\Z/d)^{d-2}$. 
Since $X$ has a $k$-point, the map $\Br(k)\to\Br(X)$ is injective and 
$\Br_1(X)/\Br(k)\simeq(\Z/d)^{d-2}$. 

Consider the following $d-2$ rational functions on $X$:
$$f_n:=\frac{x_0-\zeta^{2n+1}x_1}{x_0-\zeta x_1}, \quad n=1,\ldots,d-2,$$
where $\zeta$ is, as before, a primitive $2d$-th root of unity.
Let $k'/k$ be the splitting field of the polynomial
$x^d+s\in k[x]$ and let $G=\Gal(k'/k)\simeq\Z/d$. 
We claim that the cyclic algebras $\sA_n:=(k'/k,f_n)=(s,f_n)_d$, $n=1,\ldots,d-2$, are Azumaya
algebras on $X$ whose classes generate $\Br_1(X)/\Br(k)$. This generalises \cite[Thm.~1.1]{U}.
Moreover, if there is a subfield $k_0\subset k$ such that $k=k_0(t)$ 
has transcendence degree 1 over $k_0$, we have $\Br(X)=\Br_1(X)$ by Theorem \ref{n},
so in this case $\Br(X)/\Br(k)\simeq(\Z/d)^{d-2}$.

To prove the claim consider the rational map that projects $X$ to $\P^1_k$ with coordinates 
$(x_0:x_1)$.
It gives rise to a dominant morphism $\varphi\colon\tilde X\to\P^1_k$, where $\tilde X$ is the
blowing-up of $X$ in the subset $x_0=x_1=0$. We have $\Br(\tilde X)\cong\Br(X)$. The singular
fibres of $\varphi$ are above the $k$-points $(\zeta^{2n+1}:1)$, for $n=0,\ldots,d-1$;
the geometrically irreducible components of
these fibres are the strict transforms of the $d^2$ lines on $X_{\bar k}$ given by
$x_0^d+x_1^d=x_2^d+sx_3^d=0$.
Each singular fibre splits over $k'$ into a union of $d$ lines with a simply
transitive action of $G$. 

Recall that the {\em vertical Brauer group} associated to $\varphi$ is defined as 
$$\Br_\mathrm{vert}(\tilde X)=\Ker[\Br(\tilde X)\to\Br(\tilde X_\eta)/\varphi^*\Br(k(\eta))],$$
where $k(\eta)=k(\P^1_k)$ and $\tilde X_\eta$ is the generic fibre of $\varphi$. 
From \cite[Prop.~11.1.5]{CTS21} we see that $\Br_\mathrm{vert}(\tilde X)/\Br_0(\tilde X)$
consists of the pullbacks of the elements of $\Br(k(\eta))$ ramified at most at the $k$-points
$(\zeta^{2n+1}:1)$, for $n=0,\ldots,d-1$, such that the residue at each of these points is contained in
$$\Ker[\H^1(k,\Q/\Z)\to \H^1(k',\Q/\Z)]\cong \Hom(G,\Q/\Z).$$
Thus $\Br_\mathrm{vert}(\tilde X)/\Br_0(\tilde X)$ is generated by $\sA_1,\ldots,\sA_{d-1}$. 
The equation of $X$ shows that 
$\sA_1+\ldots+\sA_{d-1}=(s,t)_d\in\Br(k)$, so $\Br_\mathrm{vert}(\tilde X)/\Br_0(\tilde X)$ is generated
by $\sA_1,\ldots,\sA_{d-2}$. 

\ble \label{6july}
Let $Z$ be a smooth, proper, geometrically integral variety over a field $k$ of characteristic zero.
Let $\phi\colon Z\to \P^1_k$ be a surjective morphism with geometrically integral generic fibre. 
Let $U'\subset \P^1_k$ be an open subset, and let $U=\phi^{-1}(U')$. 
Let $P_U:=\Ker[\Pic(Z_{\bar k})\to \Pic(U_{\bar k})]$.
If the natural map $\H^1(k,P_U)\to\H^1(k,\Pic(Z_{\bar k}))$
is surjective, then $\Br_\mathrm{vert}(Z)\to\Br_1(Z)$ is an isomorphism.
\ele
\begin{proof} The spectral sequence $\H^p(k,\H^q(Z_{\bar k},\G_m))\Rightarrow\H^{p+q}(Z,\G_m)$
and a similar sequence for $U$ give rise to a commutative diagram with exact rows and columns
$$\xymatrix{&&\Br_1(U')\ar[d]\\
0\ar[r]&\Br_1(Z)\ar[r]\ar[d]&\Br_1(U)\ar[d]\\
\H^1(k,P_U)\ar[r]&\H^1(k,\Pic(Z_{\bar k}))\ar[r]&\H^1(k,\Pic(U_{\bar k}))}$$
Here we used that $\Br_1(U')\cong\H^2(k,\bar k[U']^\times)$ and that 
$\phi^*\colon\bar k[U']\to \bar k[U]$ is an isomorphism of $\Ga_k$-modules 
since the generic fibre of $\phi$ is proper and geometrically integral, and $\phi(U)=U'$.
The lemma follows from the diagram.
\end{proof}

\medskip

We can apply Lemma \ref{6july} to $Z=\tilde X$, $\phi=\varphi$, and $U'\subset\P^1_k$ given by
$x_0(x_0^d+x_1^d)\neq 0$.
Indeed, the morphism $\tilde X\to X$ induces a map $\tau\colon P_U\to \Lambda$, so that
we have a commutative diagram
$$\xymatrix{\H^1(k,P_U)\ar[r]^{\tau_*}\ar[d]&\H^1(k,\Lambda)\ar[d]\\
\H^1(k,\Pic(\tilde X_{\bar k}))\ar[r]^\sim&\H^1(k,\Pic(X_{\bar k}))}$$
The right hand vertical map is an isomorphism by Theorem \ref{Lambda}, so it suffices
to prove that $\tau_*$ is surjective. 
In the bottom 4-term sequence of (\ref{duo}) we can quotient
the second and third terms by the trivial $G$-module 
$\bar k[U]^\times/\bar k^\times\cong\bar k[U']^\times/\bar k^\times\cong\Z[V_f]\simeq\Z^d$. 
The second term becomes $\Z[V_g]=\Z[G]$.
The third term $\Z[G]^d\oplus\Z$ is identified with the group of divisors on 
$\tilde X_{\bar k}$ supported in $\tilde X_{\bar k}\setminus U_{\bar k}$, so its quotient is $P_U$.
Thus we obtain an exact sequence of $G$-modules
$$0\to \Z[G]/\Z\to P_U\xrightarrow{\tau} \Lambda\to 0.$$
From this we obtain an exact sequence
$$\H^1(G,P_U)\xrightarrow{\tau_*}\H^1(G,\Lambda)\to \H^2(G,\Z[G]/\Z)=0.$$
The inflation maps $\H^1(G,P_U)\to \H^1(k,P_U)$ and $\H^1(G,\Lambda)\to \H^1(k,\Lambda)$
are isomorphisms, so this proves the claim.

\brem
{When $d=3$ and $4$, the above examples may be compared to previous results of Colliot-Thélène--Kanevsky--Sansuc and Bright. In the cubic case, Examples 1 and 2 agree with the first case and Example 3 agrees with the last case in \cite[Prop.~1]{CTKS}. In the quartic case, 
Examples 1, 2, 3 agree with entries E10, E8, E20 in \cite[Table~A.5]{Bri02}, respectively.}
\erem 

\bpr \label{suppl}
Let $k$ be a field containing $\Q(\mu_{2d})$. Assume that $s,t\in k^\times$ satisfy the following conditions:

\textrm{(i)} the subgroup of $k_L^\times/(k_L^\times)^d$
generated by the classes of $s$ and $t$ is isomorphic to $(\Z/d)^2$;

\textrm{(ii)} the class of the cyclic algebra $(s,t)_d$ in $\Br(k)$ has order $d$.

\noindent Let $X\subset\P^3_k$ be the surface
\[
x_0^d-sx_1^d-tx_2^d+stx_3^d=0.
\]
Then the kernel of $\Br(k)\to\Br(X)$ is generated by the class of $(s,t)_d$.
For example, this holds when $k=k_0(s,t)$ is the purely transcendental extension 
of degree $2$ of a subfield $k_0\subset k$.
\epr
\bprf Let us prove that the class of $(s,t)_d\in\Br(k)$ goes to zero in $\Br(X)$. 
Note the kernel of $\Br(k)\to\Br(X)$ is equal to the kernel of $\Br(k)\to\Br(U)$,
where $U\subset X$ is the open subscheme given by $(x_0^d-sx_1^d)x_0\neq 0$.
Write $K=k(\sqrt[d]{s})$ and let $V$ be the $k$-torsor for the norm 1 torus
$\mathrm{R}^1_{K/k}(\G_m)$ given by $\mathrm{N}_{K/k}(\mathbf{x})=t$.
There is an obvious morphism $q\colon U\to V$, hence $\Br(k)\to\Br(U)$
is a composition $\Br(k)\to \Br(V)\to\Br(U)$. Thus it is enough to show that 
the class of $(s,t)_d$ goes to zero in $\Br(V)$. By \cite[Prop.~7.1.12]{CTS21}
the variety $V$ is birationally equivalent to the Severi--Brauer variety attached to
$(s,t)_d$. By a theorem of Amitsur, see \cite[Thm.~5.4.1]{GSz}, 
the class of $(s,t)_d$ generates the kernel of $\Br(k)\to \Br(V)$.

In view of the spectral sequence (\ref{ss}) and condition (ii)
it remains to show that the index of $\Pic(X)$ in $\Pic(X_{\bar k})^{\Ga_k}$ divides $d$.
By condition (i) and $(\mathrm{b}_X)$ in the proof of Theorem \ref{Lambda} we have 
$\Pic(X_{\bar k})^{\Ga_k}=\Lambda^{\Ga_k}$.
The Galois group $\Ga_k$ acts on $\Div_{Z_{\bar k}}(X_{\bar k})$, and hence
also on $\Lambda$, via its quotient $\Z/d$, hence $\Lambda^{\Ga_k}=\Lambda^{\Z/d}$.
From the first exact sequence of (\ref{dve}) we obtain $\H^1(\Z/d,M)\simeq\Z/d$.
Now the second exact sequence of (\ref{dve}) gives an exact sequence
\[
\Div_{Z_{\bar k}}(X_{\bar k})^{\Ga_k}\to \Lambda^{\Ga_k}\to \Z/d\to 0.
\]
Recall that the first arrow sends a divisor to its class in $\Pic(X_{\bar k})$.
It is clear that $\Div_{Z_{\bar k}}(X_{\bar k})^{\Ga_k}=\Div_{Z}(X)$.
The map $\Div_{Z}(X)\to \Pic(X_{\bar k})^{\Ga_k}$ factors through
$\Pic(X)\to\Pic(X_{\bar k})^{\Ga_k}$, so we are done.
\eprf

\section{Fields of definition}\label{Pic}

We denote the Jacobian of a smooth projective curve $C$ by $J(C)$. 
For a separable polynomial $f(t)\in k[t]$ of degree $d$ we write $C_f$ for the plane curve over $k$ defined by $y^d=f(x_1/x_0)x_0^d$.

Recall that $L_{f,g}$ is the compositum of $k(\mu_d)$ and 
the field of definition of $\Pic(Y_{\bar k})$, where $Y$ is given by (\ref{X}) with $a=1$.
We now give a computable `upper bound' for $L_{f,g}$. The number 3 in
the following proposition can be replaced by any larger integer.

\bpr\label{Lfg}
Let $K_{f,g}$ be the compositum of $k(\mu_d)$, $k_f$, $k_g$,
and the fields of definition of the 
$3$-torsion subgroups $J(C_f)[3]$ and $J(C_g)[3]$. Then $L_{f,g}\subset K_{f,g}$.
\epr
\bprf Writing $C_g$ as the curve $v^d=g(u_1/u_0)u_0^d$, there is 
a dominant rational map from $C_f\times C_g$ to $Y$ 
sending $(x_0,x_1,y,u_0,u_1,v)$ to  $(x_0 v,x_1 v,yu_0,yu_1)$. Unraveling its structure
(see \cite[Prop.~2.4, Rem.~2.7]{KSh} which develops \cite{Sas}; note the opposite sign convention for the Tate twist in \cite{KSh}), one gets an embedding of $\Ga_k$-modules:
\[
\H^2_\et(Y_{\bar k},\Q_\ell(1))\hookrightarrow
\H^2_\et(C_{f,\bar k}\times C_{g,\bar k},\Q_\ell(1))\oplus(\Q_\ell[V_f]\otimes\Q_\ell[V_g]).
\]
It is compatible with the cycle class map, so we obtain an embedding of $\Ga_k$-modules
\[
\Pic(Y_{\bar k})\otimes\Q_\ell\hookrightarrow 
\big(\NS(C_{f,\bar k}\times C_{g,\bar k})\otimes\Q_\ell\big)
\oplus(\Q_\ell[V_f]\otimes\Q_\ell[V_g]).
\]
It is well known (see, e.g.~\cite[Prop.~5.7.3]{CTS21}) that the $\Ga_k$-module
$\NS(C_{f,\bar k}\times C_{g,\bar k})$ is isomorphic to 
$\Z^2\oplus\Hom(J(C_{f,\bar k}),J(C_{g,\bar k}))$.
Thus $L_{f,g}$ is contained in the compositum of $k(\mu_d)$, $k_f$, $k_g$, and the field
of definition of $\Hom(J(C_{f,\bar k}),J(C_{g,\bar k}))$. The latter
is contained in the field of definition of 
$\Hom(J(C_f)[n],J(C_g)[n])$, for any $n\geq 3$, see \cite[Prop.~2.3]{Sil}.
\eprf

\brem{\label{defi}
For $a\in k^\times$ let $L_{f,g,a}$ be the compositum of $k(\mu_d)$ and 
the field of definition of $\Pic(X_{\bar k})$, where $X$ is given by (\ref{X}) with
the given value of $a$. The proof of Proposition \ref{Lfg} shows that 
if $\Hom(J(C_{f,\bar k}),J(C_{g,\bar k}))=0$, then $L_{f,g,a}$
is contained in the compositum of $k(\mu_d)$, $k_f$, $k_g$.
}\erem

We now focus on $L$, the field of definition of the geometric Picard group of 
the Fermat surface $F$ of degree $d$ over $\Q$.

From \cite[\S2]{GS22}, we recall a few facts about the Galois action on $\Pic(F_{\ov\Q})$. Set $E=\Q(\mu_d)$ and fix a primitive $d$-th root of unity $\e\in E$.
The group $G\simeq(\Z/d)^3$ generated by $u_1,u_2,u_3$ acts by automorphisms on $F_{\ov\Q}$, 
so that $u_i$ multiplies $x_i$ by $\e$.
Define $S_\flat$ as the set of quadruples $(l_0,l_1,l_2,l_3)$ in $\{1,\dots,d-1\}^4$ such that $\sum_{i=0}^3 (tl_i)\bmod d=2d$ for all $t\in(\Z/d)^\times$, where $(tl_i)\bmod d$
is understood as an integer in $[1,d-1]$.
An element $(l_i)\in S_\flat$ gives rise to a character $\chi\colon G\to \mu_d$ such that
$\chi(u_i)=\e^{l_i}$.
One then has a decomposition into $1$-dimensional $G$-modules
\[
\Pic(F_{\ov\Q})\otimes_\Z E=Eh\oplus\bigoplus_{\chi\in S_\flat} V_\chi,
\]
where $V_\chi\subset\Pic(F_{\ov\Q})\otimes_\Z E$ is the subspace on which $G$ acts as $g(v)=\chi(g)v$ for any $g\in G$ and $v\in V_\chi$ (cf.\ Step 5 in the proof of Theorem~\ref{Lambda}). Let us call $\chi\in S_\flat$ \emph{primitive} if it satisfies $\gcd(l_i)=1$. 

For $(l_i)\in S_\flat$ set 
$m=\gcd(l_i)$, $d'=d/m$, and $\chi'=(l_0',l_1',l_2',l_3')$ where $l_i'=l_i/m$. 
We have a morphism $f\colon F\to F'$, where $F'$ is the Fermat surface of degree $d'$,
given by $f(x_0,x_1,x_2,x_3)=(x_0^m,x_1^m,x_2^m,x_3^m)$. 
Then $\chi'$ is primitive and 
$f^*\colon V_{\chi'}\tilde\lra V_{\chi}$ is an isomorphism of $\Ga_k$-modules.

Let $\mathfrak p$ be a prime of $E$ which does not divide $d$. Let $\psi\colon\F_{\mathfrak p}^\times\to\mu_d$ be the $d$-th power residue symbol. This is the multiplicative character which assigns to $x\in\F_{\mathfrak p}^\times$ the unique $d$-th root of unity in $E$ whose reduction in $\F_{\mathfrak p}$ equals $x^{(N({\mathfrak p})-1)/d}$, 
where $N({\mathfrak p})=|\F_{\mathfrak p}|$. For $\chi\in S_\flat$, define the \emph{Jacobi sum}
\[
J_{\chi}(\mathfrak p)=
\sum_{x_1+x_2+x_3=0}\psi(x_1)^{l_1}\psi(x_2)^{l_2}\psi(x_3)^{l_3}.
\]
Then $\Frob_\mathfrak p\in\Gal(\ov\Q/E)$ acts on $V_\chi$ via multiplication by\footnote{In 
\cite[\S 2.5]{GS22} $\psi(-1)$ should be replaced by $\psi(-1)^r$ on p.~1320, lines 8, 14,
in Proposition 2.9, and on p.~1321, line $-7$.}
 
\[
h_{\chi}(\mathfrak p)=\psi(-1)^{l_0}N(\mathfrak p)^{-1} J_\chi(\mathfrak p).
\]
For an ideal 
$\mathfrak m$ of $E$, let $I(\mathfrak m)$ be the group of fractional ideals coprime to $\mathfrak m$. 
One can multiplicatively extend the domain of $h_\chi$ to $I(d\mathfrak m)$.
By a \emph{Größencharacter} modulo $\mathfrak m$ in the sense of \cite[Ch.~VII, Def.~6.1]{N}, we understand a character
\[
h:I(\mathfrak m)\to S^1=\left\{z\in\C\middle|\,|z|=1\right\}
\]
such that there exist characters
\[
h^\mathrm{fin}:(O_E/\mathfrak m)^\times\to S^1, \quad\quad h^\infty:(E\otimes_\Q \R)^\times\to S^1
\]
satisfying $h((a))=h^\mathrm{fin}(a)h^\infty(a)$ for all $a\in O_E$ coprime to $\mathfrak m$.
(Such Grö\-ßen\-charac\-ters bijectively correspond to unitary characters of the id\`ele class group of $E$.) For a given Größencharacter, the smallest $\mathfrak m$ as above is called its \emph{conductor}.

By a \emph{generalised Dirichlet character} modulo $\mathfrak m$, we understand a character of the ray class group $\Cl(\mathfrak m)=I(\mathfrak m)/P(\mathfrak m)$, where $P(\mathfrak m)$
is the group of principal ideals which are coprime to $\mathfrak m$.

\bpr\label{field}
 Let $F$ be the Fermat surface of degree $d$ over $\Q$. Let $L$ be the field of definition of 
$\Pic(F_{\ov \Q})$.

 \textrm{(i)}  If $d\leq 4$ or $\mathrm{gcd}(6,d)=1$, then $L=\Q(\mu_{2d})$.
 
 \textrm{(ii)} If $d$ is not divisible by any of the integers between $12$ and $180$ in the table of 
\textrm{\cite[p.~727]{Sh81}}, then $L=\Q(\mu_{2d},2^{2/d})$ if $d$ is $2$ or $4$ modulo $6$,
$L=\Q(\mu_{2d},3^{3/d})$ if $d$ is $3$ modulo $6$, and 
$L=\Q(\mu_{2d},2^{2/d}, 3^{3/d})$ if $d$ is divisible by $6$.

\epr
\bprf
\textrm{(i)} Any line
in $F_{\ov\Q}$ is a unique effective divisor in its class because it has negative self-intersection 
$2-d$ by the Riemann--Roch theorem for surfaces.
Thus from the equations of the $3d^2$ obvious lines in $F_{\ov\Q}$ we infer that
$\Q(\mu_{2d})\subset L$. By Shioda \cite[Thm.~7]{Sh81},  $\Pic(F_{\ov\Q})\otimes_\Z\Q$ is generated by the classes of these lines if $d\leq 4$ or $\mathrm{gcd}(6,d)=1$.

\textrm{(ii)} Let $S_\mathrm{ind}\subset S_\flat$ (for \emph{indecomposable}) be the subset of quadruples 
$(l_0,l_1,l_2,l_3)$ such that $l_i+l_j\neq d$ for all $i\neq j$.
Let $\Lambda'\subseteq\Pic(F_{\ov \Q})$ be the subgroup generated by the $3d^2$ lines. 
As explained in \cite[Proof of Thm.~7]{Sh81}, the classes of these lines span $\bigoplus_{\chi\in S_\flat\setminus S_\mathrm{ind}}V_\chi$. Thus
\[
\Pic(F_{\ov \Q})\otimes_\Z E=(\Lambda'\otimes_\Z E)\oplus\bigoplus_{\chi\in S_\mathrm{ind}}V_\chi.
\]
Following \cite[Lem.~1]{Sh81}, let $S_\mathrm{reg}\subseteq S_\mathrm{ind}$ (for \emph{regular}) be the subset of quadruples which, possibly after a permutation, can be written in one of the following forms:

{(a)} if $d=2d'$ is even,
\begin{align*}
 (i,d'+i,d-2i,d'),\ & 1\leq i\leq d'-1, \ i\neq d'/2, \ \gcd(i,d')=1;\\
 (i,d'+i,d'+2i,d-4i),\ & 1\leq i\leq (d'-1)/2,\  i\neq d'/3, \  \gcd(i,d')=1;\\
 (i,d'+i,2i-d',2d-4i),\ & (d'+1)/2\leq i\leq d'-1, \ i\neq d/3, \ \gcd(i,d')=1;
\end{align*}

{(b)} if $d=3d''$ is divisible by $3$,
\begin{align*}
 (i,d''+i,2d''+i,d-3i),\ & 1\leq i\leq d''-1, \  i\neq d''/2, \ \gcd(i,d'')=1.
\end{align*}

In \cite[Thm.~1--3]{ASh}, Aoki and Shioda find for each $\chi\in S_\mathrm{reg}$ an explicit equation describing a divisor whose class generates $V_{\chi}$\footnote{\cite{ASh} writes the elements of $S_\mathrm{reg}$ in a slightly different but equivalent way.}. These divisors are defined over $\Q(\mu_{2d},2^{2/d})$ (for $d$ divisible by $2$) or $\Q(\mu_{2d},3^{3/d})$ (for $d$ divisible by $3$) and have negative self-intersection, hence they are unique effective divisors in their classes. 
By \cite[Thm.~C]{Aok} and \cite[Eq.~$(9)$]{Sh81}, if $\chi\in S_\mathrm{ind}$ is primitive and $d$ is not amongst the integers between $12$ and $180$ in the table of \cite[p.~727]{Sh81}, then $\chi\in S_\mathrm{reg}$.
It then follows that
if $d$ is not divisible by any of the integers between $12$ and $180$ in the table of  \cite[p.~727]{Sh81},  
then $S_\mathrm{ind}=S_\mathrm{reg}$. The claim (ii) follows.
\eprf

\bpr \label{D}
Write $w=|(O_E^\times)_\mathrm{tors}|$, so that $w=d$ if $d$ is even, and $w=2d$ if $d$ is odd. Let $\Delta$ be the subgroup of $E^\times/E^{\times w}$
generated by $O_E^\times$ and the elements $1-\e^{d/p}$ 
for all prime divisors $p|d$. Then $L$ is contained in $K=E(\Delta^{1/w})$.
\epr
\bprf Reducing the problem as in the proof of Proposition \ref{field} (ii), it remains to determine the image of the action of $\Gal(\ov\Q/E)$ on $V_\chi$ for all $\chi=(l_0,l_1,l_2,l_3)$ in $S_\mathrm{ind}\setminus S_\mathrm{reg}$.

It was first shown by Weil \cite[Theorem]{W2} that $h_\chi$ is a Größencharacter of conductor dividing $d^2$, hence dividing $c\coloneq 2w^2$. Since $\Gal(\ov\Q/E)$ acts on $V_\chi\subset\Pic(F_{\ov\Q})\otimes_\Z E$ through a finite abelian quotient, it follows by \cite[Ch.~VII, Prop.~6.9]{N} that $h_\chi$ factors through a generalised Dirichlet character $h^\mathrm{fin}_\chi$ modulo~$c$.

By class field theory, the action of $\Gal(\ov\Q/E)$ on $V_\chi$ thus factors through the Galois group $\Gal(E(c)/E)\simeq\Cl(c)$ of the ray class field $E(c)$ of conductor $c$. By Lemma~\ref{kumext} below, we have $K\subset E(c)$ and hence the Artin map $a_{K/E}:I(c)\to \Gal(K/E)$ factors through a homomorphism $a^\mathrm{fin}_{K/E}:\Cl(c)\to\Gal(K/E)$. It remains to show that
\begin{equation}\label{cft}
 \Ker(a^\mathrm{fin}_{K/E})\subset \Ker(h^\mathrm{fin}_\chi)
\end{equation}
in the finite ray class group $\Cl(c)$. For any given $\chi$, the determination of $\Ker(h^\mathrm{fin}_\chi)$ is a finite calculation with Jacobi sums which can be performed as in \cite[\S4.1]{W} with customary computer algebra packages (see Remark~\ref{magma}); and so is the determination of 
$\Ker(a^\mathrm{fin}_{K/E})$ using the characterisation of the Artin map of a Kummer extension from \cite[\S2]{Fie}.

It suffices to verify (\ref{cft}) for each primitive $\chi\in S_\mathrm{ind}\setminus S_\mathrm{reg}$ and all $d$. 
By the aforementioned \cite[Thm.~C]{Aok}, such $\chi$ may only exist for those $12\leq d\leq 180$ that are listed in the table of \cite[p.~727]{Sh81}. Checking this finite list of $\chi$ completes the proof of (iii).
\eprf

\ble\label{kumext}
We have $K\subset E(c)$.
\ele
\bprf It is enough to show that if
$x\in O_E$ divides $d\geq 3$, then the conductor of the Kummer extension 
$E(\sqrt[w]{x})$ of $E$ divides $c$.
The relative discriminant $\Delta_{E(\sqrt[{w}]{x})/E}$ divides the relative discriminant ideal of the polynomial $t^{w}-x$, which is $({w}^{w}x^{{w}-1})$. 
A formula in Hasse's lectures on class field theory \cite[Sätze (163), (164)]{Has} says that the conductor of $E(\sqrt[{w}]{x})$ divides
\[
\prod_{\mathfrak p\nmid {w}}\mathfrak p\prod_{\mathfrak p|{w}}\mathfrak p^{\left(v_p({w})+\frac{1}{p-1}\right)e(\mathfrak p/p)+1},
\]
where both products are taken over all primes $\mathfrak p\subset O_E$ dividing $\Delta_{E(\sqrt[{w}]{x})/E}$, and $e(\mathfrak p/p)$ is the ramification index of $\mathfrak p$ over the underlying prime number $p$.
The first product is vacuous. Now if $v_p(w)>1$ or $p>2$, then $e(\mathfrak p/p)/(p-1)+1\leq v_p(w)e(\mathfrak p/p)$. If $p=2$ and $v_2(w)=1$, then $\left(v_2({w})+\frac{1}{2-1}\right)e(\mathfrak p/2)+1 = 2v_2(w)e(\mathfrak p/2)+1$. Thus the conductor of $E(\sqrt[{w}]{x})$ divides 
\[\prod_{\mathfrak p|{w}}\mathfrak p^{2v_p(w)e(\mathfrak p/p)}\prod_{\mathfrak p|2}\mathfrak p,\]
hence it divides $c$.
\eprf

\brem\label{magma}{
The calculation in Proposition \ref{D}
leads to an explicit description of the exact field of definition $L$ of the geometric
Picard group of the Fermat surface of degree $d$, for all positive integers $d$. 
This is achieved by computing, 
for each $12\leq d\leq 180$ in the table of \cite[p.~727]{Sh81} and for each primitive
$\chi\in S_\flat$ in degree $d$, 
the exact kernel of $h^\mathrm{fin}_\chi$ 
and the corresponding abelian extension $L_\chi/E$. Let $M_d$ be the compositum of the fields 
$L_\chi$, for all primitive $\chi\in S_\flat$ in degree $d$.
Then $L$ is the compositum of the fields $E=\Q(\mu_{2d})$, $\Q(2^{2/d})$ if $2\mid d$, $\Q(3^{3/d})$ if $3\mid d$, and the fields $M_{d'}$ for all $d'|d$, $12\leq d'\leq 180$. The table of the fields $M_d$ and the associated Magma code \cite{Magma}, which also verifies Proposition~\ref{D}, are available at 
\href{https://github.com/dgvirtz/fermatdeffield}{github.com/dgvirtz/fermatdeffield}.
}
\erem


\noindent Department of Mathematics, University of Glasgow, University Place,
Glasgow,  G12~8QQ United Kingdom

\medskip

\noindent \texttt{damian.gvirtz@glasgow.ac.uk}

\bigskip

\noindent Department of Mathematics, South Kensington Campus,
Imperial College London, SW7 2BZ England, U.K. -- and --
Institute for the Information Transmission Problems,
Russian Academy of Sciences, 19 Bolshoi Karetnyi, Moscow, 127994
Russia

\medskip

\noindent \texttt{a.skorobogatov@imperial.ac.uk}


\begin{thebibliography}{99}







\bibitem{Aok}
Aoki, N.: On some arithmetic problems related to the Hodge cycles on the Fermat varieties.
Math. Ann. \textbf{266}, 23--54 (1983). Erratum: Math. Ann. \textbf{267}, 572 (1984)

\bibitem{ASh}
Aoki, N., Shioda, T.: Generators of the Néron--Severi group of a Fermat surface.
In: Arithmetic and geometry, vol.~I, pp. 1--12,
Progr. Math. \textbf{35}, Birkhäuser, Boston, MA, 1983.
Correction and supplement to the paper ``Generators of the N\'eron--Severi group of a Fermat surface'': Comment. Math. Univ. St. Pauli \textbf{59}, 65--75 (2010)

\bibitem{Bri02}
Bright, M.: Computations on diagonal quartic surfaces.
Ph.D. thesis, University of Cambridge (2002)

\bibitem{Bri11}
Bright, M.: The Brauer--Manin obstruction on a general diagonal quartic surface.
Acta Arith. \textbf{147}, 291--302 (2011)

\bibitem{BBL}
Bright, M.~J., Browning, T.~D., Loughran, D.: Failures of weak approximation in families.
Compositio Math. \textbf{152}, 1435--1475 (2016)

\bibitem{CTKS}
Colliot-Thélène, J.-L., Kanevsky, D., Sansuc, J.-J.: Arithmétique des surfaces cubiques diagonales.
In: Diophantine approximation and transcendence theory, pp. 1--108,
Lecture Notes in Mathematics \textbf{1290}, Springer-Verlag, Berlin, 1987

\bibitem{CTS21}
Colliot-Thélène, J.-L., Skorobogatov, A.~N.: The Brauer--Grothendieck group.
Ergebnisse der Mathematik und ihrer Grenzgebiete, 3. Folge, Band 71,
Springer, Cham, 2021

\bibitem{Deg15}
Degtyarev, A.: Lines generate the Picard groups of certain Fermat surfaces.
J. Number Theory \textbf{147}, 454--477 (2015)

\bibitem{Fie}
Fieker, C.: Computing class fields via the Artin map.
Math. Comp. \textbf{70}, 1293--1303 (2001)

\bibitem{GSz}
Gille, P., Szamuely, T.: Central simple algebras and Galois cohomology. 2nd ed.
Cambridge Stud. Adv. Math. \textbf{165}, Cambridge University Press, Cambridge (2017)

\bibitem{GS22}
Gvirtz, D., Skorobogatov, A.~N.: Cohomology and the Brauer groups of diagonal surfaces.
Duke Math. J. \textbf{171}, 1299--1347 (2022)

\bibitem{Has}
Hasse, H.: Vorlesungen über Klassenkörpertheorie.
Thesaurus Math. \textbf{6}, Physica-Verlag, Würzburg (1967)

\bibitem{KSh}
Katsura, T., Shioda, T.: On Fermat varieties.
Tôhoku Math. J. \textbf{31}, 97--115 (1979)

\bibitem{Magma}
Bosma, W., Cannon, J., Playoust, C.: The Magma algebra system. I The user language.
J. Symbolic Comput. \textbf{24}, 235--265 (1997)

\bibitem{N}
Neukirch, J.: Algebraic number theory.
Grundlehren Math. Wiss. \textbf{322}, Springer-Verlag, Berlin (1999)

\bibitem{PT}
Peyre, E., Tschinkel, Yu.: Tamagawa numbers of diagonal cubic surfaces, numerical evidence.
Math. Comp. \textbf{70}, 367--387 (2001)

\bibitem{San}
Santens, T.: Diagonal quartic surfaces with a Brauer--Manin obstruction.
Compos. Math. \textbf{159}, 659--710 (2023)

\bibitem{Sas}
Sasakura, N.: On some results on the Picard numbers of certain algebraic surfaces.
J. Math. Soc. Japan \textbf{20}, 297--321 (1968)

\bibitem{Serre}
Serre, J.-P.: Lectures on the Mordell--Weil theorem. 2nd ed.
Friedr. Vieweg \& Sohn, Braunschweig (1990)

\bibitem{SGA}
Deligne, P., Katz, N.: Groupes de monodromie en géométrie algébrique, II.
Sém. géométrie algébrique du Bois-Marie 1967--69 (SGA 7 II),
Lecture Notes in Math. \textbf{340}, Springer-Verlag, Berlin (1973)

\bibitem{Sh79}
Shioda, T.: The Hodge conjecture for Fermat varieties.
Math. Ann. \textbf{245}, 175--184 (1979)

\bibitem{Sh81}
Shioda, T.: On the Picard number of a Fermat surface.
J. Fac. Sci. Univ. Tokyo Sect. IA Math. \textbf{28}, 725--734 (1981)

\bibitem{Sil}
Silverberg, A.: Fields of definition for homomorphisms of abelian varieties.
J. Pure Appl. Algebra \textbf{77}, 253--262 (1992)

\bibitem{U}
Uematsu, T.: On the Brauer group of diagonal cubic surfaces.
Quart. J. Math. \textbf{65}, 677--701 (2014)

\bibitem{U2}
Uematsu, T.: On the Brauer group of affine diagonal quadrics.
J. Number Theory \textbf{163}, 146--158 (2016)

\bibitem{W}
Watkins, M.: Jacobi sums and Grössencharacters.
Publ. Math. Besançon Algèbre Théorie Nr. (2018), 111--122

\bibitem{Wei}
Wei, D.: Strong approximation for a toric variety.
Acta Math. Sin. (Engl. Ser.) \textbf{37}, 95--103 (2021)

\bibitem{W2}
Weil, A.: Jacobi sums as “Grössencharaktere”.
Trans. Amer. Math. Soc. \textbf{73}, 487--495 (1952)

\bibitem{Weibel}
Weibel, C.: An introduction to homological algebra.
Cambridge Stud. Adv. Math. \textbf{38}, Cambridge University Press, Cambridge (1994)


\end{thebibliography}
\end{document}